\documentclass[12pt]{article}

\usepackage{preprint-layout}
\usepackage{preprint-notation}

\usepackage[utf8]{inputenc}
\usepackage{booktabs}

\title{%
The Finite Cell Method with Least Squares Stabilized Nitsche Boundary Conditions
}
\date{\today}
\author{
Karl Larsson,\
Stefan Kollmannsberger,\ 
\\
Ernst Rank,\
Mats G. Larson
}
\date{\today}

\newcommand{\supp}{\text{supp}}
\newcommand{\Span}{\text{span}}
\newcommand{\hatA}{\widehat{A}}

\begin{document}

\maketitle

\begin{abstract}
We apply the recently developed  least squares stabilized symmetric Nitsche method for 
enforcement of Dirichlet boundary conditions to the finite cell method. The least squares stabilized Nitsche method in combination with finite cell stabilization leads to a symmetric positive definite stiffness matrix and relies only on  elementwise stabilization, which does not 
lead to additional fill in. We prove a priori error estimates and bounds on the condition numbers. 
\end{abstract}

\section{Introduction}

Unfitted finite element methods allow for the geometry of the problem domain to cut through the underlying computational mesh in an arbitrary manner and avoids the use of standard meshing tools. 
Over recent years several methods that deal with the challenges of unfitted FEM in slightly different ways have been proposed, see for instance \cite{MR2377802,MR1941489,BurClaHanLarMar2015,MR2497337,MR2458114}.

In this contribution we provide an analysis for an unfitted FEM, the finite cell method~\cite{MR2377802}, when combined with an improved technique for weakly imposing Dirichlet boundary conditions. Specifically, we address the following two central issues that naturally arise when developing unfitted FEMs:

\paragraph{1. Enforcing Dirichlet Boundary Conditions.}

In contrast to classical FEMs the degrees of freedom in unfitted FEMs are no longer nodal values coinciding with the boundary. A very elegant solution for dealing with this issue is weak enforcement via Nitsche's method, first introduced in 1971 by Nitsche in his seminal work \cite{Nit1971}. Weak enforcement of Dirichlet boundary conditions does seem to be a perfect match with unfitted methods and the elegance of this particular method shows in optimal order error bounds using only a $h^{-1}$ scaling of the penalty parameter.
However, the stability analysis in \cite{Nit1971} utilizes an inverse inequality that does not hold by default for unfitted finite elements, whereby additional stabilization is required for a complete analysis.

\paragraph{2. Ensuring Stability and Accuracy.}
While there are clear modeling benefits of an unfitted FEM approach, a drawback is that there may exists basis functions in our computational grid whose support has a negligible intersection with the domain. This takes its expression in the form of two stability issues that are related but not the same:
\begin{itemize}
\item
\emph{Stability of Formulation.} For elliptic PDE, as we consider in this contribution, this means a guaranteed coercive method. If the method is not coercive, we cannot be sure we are actually solving the correct problem, even if we may solve the resulting linear system of equations without any hassle.
\item
\emph{Stability of Linear System of Equations.}
This means there exists an upper bound on the condition number for the stiffness matrix.
\end{itemize}
To mitigate these stability issues several approaches have been proposed in literature, for instance;
ghost penalty stabilization \cite{Bur2010,BurClaHanLarMar2015};
element merging \cite{JohLar2013};
removal of problematic basis functions \cite{ElfLarLar2018};
and finite cell stabilization \cite{DauDusRan2015,MR3311666}.

\paragraph{The Finite Cell Method.}
As a framework for unfitted FEMs, addressing problems of discretization, quadrature, and stability, the finite cell method (FCM) was originally introduced in \cite{MR2377802} and has since successfully been applied to a wide variety of applications, for instance in solid mechanics \cite{MR2458114,MR2896903}, structural analysis \cite{MR3358026,MR3944486}, and fluid-structure interaction \cite{MR3310316}. 
In this work, the focus is on the stability issues and the mathematical analysis of the method.
The principal stabilization technique in FCM is based on adding a small virtual stiffness in parts of elements that are exterior to the domain, and this we refer to as \emph{finite cell stabilization}. In the case of Neumann conditions, an analysis of finite cell stabilization was provided in \cite{DauDusRan2015}, while we in the present contribution provide an analysis in the case of weakly enforced Dirichlet conditions.

\paragraph{Summary of Contribution.}
In this paper we develop a new symmetric finite cell method which utilizes the recently developed
least-squares stabilized Nitsche formulation (LS-Nitsche) \cite{ElfLarLar2019} for weakly enforcing
the Dirichlet boundary conditions.
We work in a higher regularity space $V = H^2(\Omega)$ setting and we therefore consider
$C^1$ splines constructed on a background grid as our approximation space $V_h$.
The following key considerations are taken into account:

\begin{itemize}
\item \emph{Stability of Formulation.}
The use of LS-Nitsche \cite{ElfLarLar2019} appends the standard Nitsche formulation by certain least squares terms in the vicinity of the boundary as well as control over the tangent gradient along the boundary, yielding a formulation which is stable not only on the discrete space $V_h$ but on $V = H^2(\Omega)$. Effectively, this means the method is coercive regardless of the cut situation even without any finite cell stabilization.

\item \emph{Stability of Linear System of Equations.}
Since coercivity is ensured via LS-Nitsche, the focus of the finite cell stabilization, i.e. the virtual stiffness added inside cut elements outside of the domain, is to guarantee stability of the linear system of equations.
If finite cell stabilization is not included, the norm in which coercivity is proven will only be a seminorm on $V_h$ since $V_h$ in certain cut situations may include functions whose support has an arbitrarily small intersection with the domain. This, in effect, gives a stiffness matrix which is only guaranteed to be symmetric positive semidefinite, while including finite cell stabilization will move the stiffness matrix to be symmetric positive definite (SPD).
We include an error analysis, that identifies the suitable choice of the finite cell stabilization parameter. Note that thanks to LS-Nitsche we already have coercivity and only very weak virtual stiffness need to be added to ensure that the stiffness matrix is SPD. Our results complements 
\cite{DauDusRan2015} where Neumann conditions are analyzed.

\item \emph{Implementation and Use.}
In the method, all terms are assembled element-wise making for a very straightforward implementation.
For this to hold the added smoothness of the spline space is utilized in the derivations of the least-squares terms.
A feature of LS-Nitsche that simplifies the use of the method is that the choice of penalty parameter
does not depend on the cut situation. Furthermore, a modest value for the penalty parameter may be chosen which allows complex geometries even at the element level, without locking.

\end{itemize}

\paragraph{Outline.} The paper is organized as follows. In Section 2 we introduce our method; an unfitted finite element method featuring least-squares stabilized Nitsche boundary conditions and finite cell stabilization. In Section 3 we prove stability and error estimates in the energy and $L^2$ 
norms. In Section 4 we derive a bound on the condition number. In Section 5 we present illustrating numerical examples. Finally, in Section 6 we present our conclusions.

\section{The Finite Cell Method with LS-Nitsche}

\subsection{Model Problem}

\begin{figure}\centering
\begin{subfigure}[t]{0.40\linewidth}\centering
\includegraphics[width=0.85\linewidth]{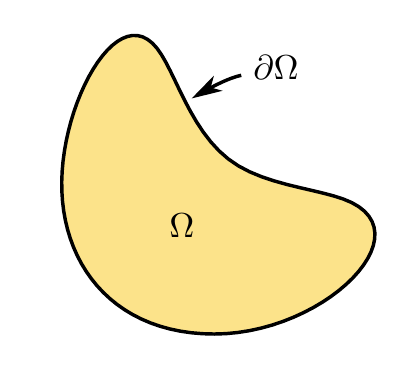}
\subcaption{Problem domain}\label{bean:a}
\end{subfigure}
\begin{subfigure}[t]{0.40\linewidth}\centering
\includegraphics[width=0.85\linewidth]{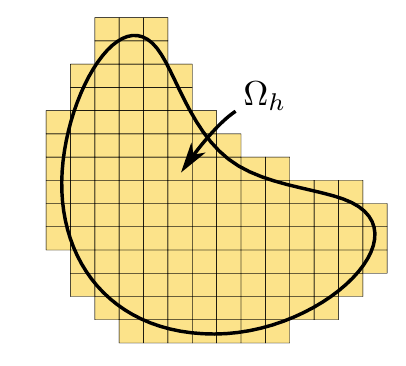}
\subcaption{Domain of active elements}\label{fig:bean:b}
\end{subfigure}
\caption{
(a) The domain $\Omega$ of the model problem.
(b) The domain of all active elements $\Omega_h$, i.e. the union of all elements with a non-empty intersection with the domain $\Omega$.
}
\label{fig:bean0}
\end{figure}

Let $\Omega$ be a domain in $\IR^d$ with smooth boundary $\partial \Omega$, such as illustrated in Figure~\ref{bean:a}, 
and consider the problem: find $u: \Omega \rightarrow \IR$ such that 
\begin{alignat}{3}\label{eq:problem-strong-a}
-\Delta u &=f & \qquad &\text{in $\Omega$}
\\ \label{eq:problem-strong-b}
u &=g & \qquad &\text{on  $\partial \Omega$}
\end{alignat}

For sufficiently regular data there exists a unique solution to this problem and we will be 
interested in higher order methods and therefore we will always assume that the solution 
satisfies the regularity estimate 
 \begin{equation}\label{eq:regularity}
 \| u \|_{H^{s+2}(\Omega)} \lesssim 
 \|f \|_{H^{s}(\Omega)} 
 + 
 \|g \|_{H^{s-1/2}(\partial\Omega)}
 \end{equation}
 for some $s \geq 0$.
Here and below $a \lesssim b$ means that there is 
a positive constant $C$ independent of the mesh parameter $h$ and the cut situation such that $a \leq C b$.

\subsection{The B-Spline Spaces}

\begin{figure}
\centering
\includegraphics[width=0.9\linewidth]{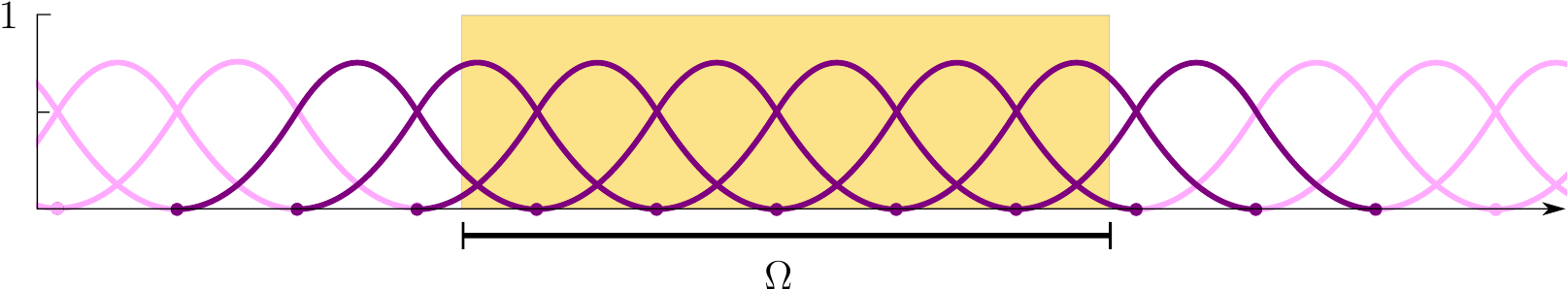}
\caption{$C^1 Q^2(\IR)$ B-spline basis functions in one dimension. The set $B$ of basis functions with non-empty support in $\Omega$ is indicated in deep purple. Note that basis functions crossing the boundary 
of $\Omega$ are defined analogously to interior basis functions.}
\label{fig:splines}
\end{figure}
\paragraph{Definitions.}
\begin{itemize}
\item Let $\widetilde{\mcT}_{h}$, $h \in (0,h_0]$, be a family of uniform tensor 
product meshes in $\IR^d$ with mesh parameter $h$. 
\item Let $\widetilde{V}_{h}=C^{p-1}Q^p(\IR^d)$ be the space of 
$C^{p-1}$ tensor product B-splines of order $p$ defined on 
$\widetilde{\mcT}_{h}$.
Such functions are easily generated using the Cox-de Boor recursion formula.
Let $\widetilde{B} = \{\varphi_i\}_{ i\in \widetilde{I}}$ be 
the standard basis in $\widetilde{V}_h$, where $\widetilde{I}$ is an index set.
\item Let $B = \{ \varphi \in \widetilde{B} \, : \,
\supp(\varphi) \cap \Omega \neq \emptyset \}$ be the set of basis functions with 
support that intersects $\Omega$. Let $I$ be an index set for $B$. Let $V_h = \Span\{B\}$, 
$\mcT_h = \{ T \in \widetilde{\mcT}_h : T \subset \cup_{\varphi \in B } \supp(\varphi) \}$, and 
$\Omega_h = \cup_{T \in \mcT_h} T$. An illustration of the basis functions in one spatial dimension is given in Figure~\ref{fig:splines} and $\Omega_h$ is illustrated in Figure~\ref{fig:bean:b}.
\item We will only consider $p \geq 2$ corresponding to at least $C^1$ splines. We 
then have $V_h|_\Omega \subset V = H^2(\Omega)$. 

\end{itemize}

\subsection{The  Method}

\begin{tcolorbox}[colback=white,boxrule=1.5pt]
\paragraph{Method.}
Find 
$u_{h,\alpha} \in V_{h}$ such that 
\begin{equation}\label{eq:method}
A_{h,\alpha}(u_{h,\alpha},v) = L_h(v)\qquad \forall v \in V_{h}
\end{equation}
\tcblower
The forms are defined as follows:
\begin{align}\label{eq:Ah}
A_{h,\alpha}(v,w) &= a_h(v,w) +
\overbrace{\alpha(\nabla v, \nabla w)_{\Omega_h \setminus \Omega}}^{\mathclap{\text{Finite cell stabilization}}}
\\ 
&\qquad
\underbrace{
- (\nabla_n v,w)_{\partial \Omega} 
- (v, \nabla_n w)_{\partial \Omega} 
+\beta b_h(v,w)
}_{\mathclap{\text{Nitsche BC terms}}}
\\ \label{eq:ah}
a_h(v,w)&= (\nabla v, \nabla w)_\Omega
+
\underbrace{
\tau h^2(\Delta v,\Delta w)_{\mcT_{h,\partial\Omega} \cap \Omega}
}_{\mathclap{\text{New LS interior term}}}
\\ \label{eq:bh}
b_h(v,w)&=
  (2 + \tau^{-1} ) h^{-1} (v,w)_{\partial \Omega}
+
\underbrace{
2 h (\nabla_T  v, \nabla_T w)_{\partial \Omega} 
}_{\mathclap{\text{New LS boundary term}}}
\\ \label{eq:Lh}
L_h(v) &= (f, v)_\Omega
\underbrace{
- \tau h^2 (f,\Delta v)_{\mcT_{h,\partial\Omega}\cap \Omega}
}_{\mathclap{\text{New LS interior term}}}
\underbrace{
- (g,\nabla_n v)_{\partial\Omega} + \beta b_h(g,v)
}_{\mathclap{\text{Nitsche BC terms}}}
\end{align}


\end{tcolorbox}

\paragraph{Parameters.}
The included parameters are:
\begin{itemize}
\item
$\alpha \geq 0$ is the finite cell stabilization 
parameter, which we will define as $\alpha \sim h^{2p -1}$ to obtain a scheme with optimal order convergence. This stabilization provides a small virtual stiffness in exterior parts of cut elements which, for the present method, is crucial for guaranteeing stability of the linear system of equations.

\item  $\beta>0$ is a penalty parameter of moderate size, for instance $\beta = 5$ is sufficient, see \cite{ElfLarLar2019}.

\item $\tau > 0$ can be chosen in order to optimize accuracy, where typical values are in the range 0.01--1, see \cite{ElfLarLar2019}. The effective 
parameter in front of the standard Nitsche penalty term  $h^{-1}(v,w)_{\partial \Omega}$ is 
$\beta(2 + \tau^{-1})$ and thus $\tau$ trades control between the interior least squares 
term $h^2(\Delta v, \Delta w )_{\mcT_{h,\partial\Omega} \cap \Omega}$ and the standard Nitsche 
penalty term.
\end{itemize}

\begin{figure}\centering
\begin{subfigure}[t]{0.40\linewidth}\centering
\includegraphics[width=0.97\linewidth]{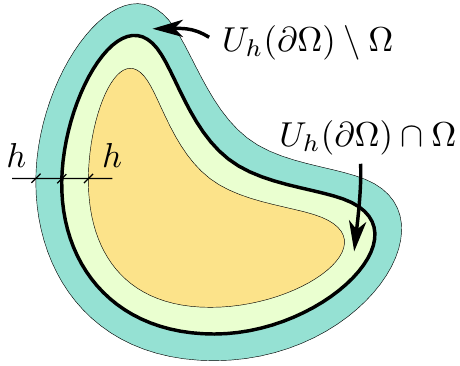}
\subcaption{}
\end{subfigure}
\begin{subfigure}[t]{0.40\linewidth}\centering
\includegraphics[width=0.85\linewidth]{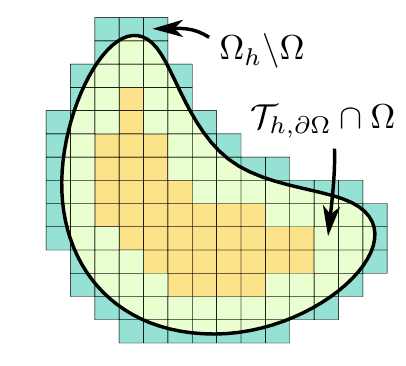}
\subcaption{}
\end{subfigure}
\caption{
Illustrations of geometrical regions used in the method and throughout the analysis. By construction $(U_h(\partial\Omega)\cap\Omega) \subseteq (\mcT_{h,\partial\Omega}\cap \Omega)$ and $(\Omega_h\setminus\Omega) \subseteq (U_h(\partial\Omega) \setminus\Omega)$.
}
\label{fig:geom-regions}
\end{figure}

\paragraph{Notation.}
We employed the following notation:
\begin{itemize}
\item $\nabla_n = n \cdot \nabla$ is the normal derivative on the boundary $\partial\Omega$ where $n$ is the exterior unit normal to $\partial\Omega$.

\item  $\nabla_T = P\nabla$ is the so-called tangent gradient on  $\partial \Omega$ where $P = (I - n\otimes n)$ is the projection onto the tangent plane of $\partial \Omega$.
\item  $\mcT_{h,\partial\Omega} \cap \Omega$ is the interior region of least squares stabilization, see Figure~\ref{fig:geom-regions},
where $\mcT_{h,\partial\Omega} \subset \mcT_h$ is defined
\begin{equation}\label{def:Th}
\mcT_{h,\partial\Omega} = \mcT_h(U_h (\partial \Omega))
= \{T \in \mcT_h : T \cap U_{h} (\partial \Omega) \neq \emptyset \}  
\end{equation}
and $U_h (\partial \Omega)$ is the tubular neighborhood
\begin{equation}\label{eq:Udelta}
U_h (\partial \Omega) = \bigcup_{x\in \partial \Omega} B_h(x)
\end{equation}
with $B_h(x)$ the open ball with center $x$ and radius $h$.
Integrals over $\mcT_{h,\partial\Omega} \cap \Omega$ are evaluated elementwise such that
\begin{equation}
(v,w)_{\mcT_{h,\partial\Omega} \cap \Omega} 
= \sum_{\mathclap{T \in \mcT_{h,\partial\Omega}}} (v,w)_{T\cap \Omega}
\end{equation}
\end{itemize}
\begin{rem}
In practice, $\mcT_{h,\partial\Omega}$ may be taken as the set of all elements that intersect 
the Dirichlet boundary $\partial \Omega$ and their neighbors, i.e., $\mcT_{h,\partial\Omega} = 
\mcN_h(\mcT_h(\partial \Omega))$.
\end{rem}

\paragraph{Commentary on LS Terms.}

The added least-squares terms, $h^2(\Delta v,\Delta w)_{\mcT_{h,\partial\Omega} \cap \Omega}$ and $h (\nabla_T v, \nabla_T w)_{\partial \Omega}$, 
 provide additional control on the elements in the vicinity of the boundary and along the boundary. The additional control along the boundary may be interpreted as weak 
 enforcement of  Dirichlet boundary conditions in $H^{1/2}(\partial \Omega)$, the trace space 
 of $H^1(\Omega)$, and we note that 
\begin{equation}
h^{-1} \| v \|^2_{\partial \Omega} + h \| \nabla_T v \|^2_{\partial \Omega} 
\end{equation}
is a discrete version of $\| v \|^2_{H^{1/2}(\partial \Omega)}$, which is more precise compared 
to the standard Nitsche method due to the presence of the second term.  In the case when we 
do not have any cut elements we may employ an inverse inequality to conclude that 
\begin{equation}
h^{-1} \| v \|^2_{\partial \Omega} + h \| \nabla_T v \|^2_{\partial \Omega}  \lesssim h^{-1} \| v \|^2_{\partial \Omega}
\end{equation}
and we recover the standard Nitsche penalty term. This is however not possible in the cut case 
where the additional term plays a key role. Furthermore, for $w \in H^1(\Omega)$ we have the 
standard estimate 
\begin{align}
\left| (\nabla_n v, w)_{\partial \Omega} \right| \leq 
(\| \nabla v \|^2_\Omega + \| \Delta v \|^2_\Omega)^{1/2} \| w \|_{H^1(\Omega)}
\end{align} 
and thus for the left hand side to be well posed we need control of $ \| \Delta v \|_\Omega$ in 
addition to $\| \nabla v \|_\Omega$. We note that $h^2(\Delta v,\Delta w)_{\mcT_{h,\partial\Omega} \cap \Omega}$ and $h (\nabla_T v, \nabla_T w)_{\partial \Omega}$ provides such control close to the boundary, which turns out to be enough. Finally, if we do not have cut elements an inverse inequality gives
\begin{equation}
h^2 \| \Delta v \|^2_\Omega\lesssim \| \nabla v \|^2_\Omega
\end{equation}
and thus in that case $h^2 \| \Delta v \|^2_\Omega$ is dominated 
by $ \| \nabla v \|^2_\Omega$, which is already present in the standard variational form.

\section{Analysis of the Method}

\subsection{Extension Outside the Domain}
Throughout the analysis we need to evaluate functions in $V=H^2(\Omega)$ not only within $\Omega$ but also a region exterior to $\Omega$, specifically within $U_h(\partial\Omega) \setminus\Omega$, see Figure~\ref{fig:geom-regions}.
For this we use the following known result; there is an extension operator $E:W^k_q(\Omega) \rightarrow W^k_q(\IR^d)$, $k\geq 0$ 
and $q\geq 1$, such that 
\begin{equation}\label{eq:extension-continuity}
\| E v \|_{W^k_q(\IR^d)} \lesssim \| v \|_{W^k_q(\Omega)}
\end{equation}
see \cite{Fol1995}.
In particular we use this extension on the exact solution $u$ in our derivations below, where we let $u$ on $U_h(\partial \Omega) \cup \Omega$ be defined by 
$u$ on $\Omega$ and by the extension $E u$ on $U_h(\partial\Omega) \setminus \Omega$.

\subsection{Consistency}

\begin{tcolorbox}[colback=white,boxrule=1.5pt]
\paragraph{Galerkin Orthogonality.} With $u \in V = H^2(\Omega)$ the solution to 
\eqref{eq:problem-strong-a}--\eqref{eq:problem-strong-b} and $u_{h,\alpha} \in V_h$ 
the solution to (\ref{eq:method}) it holds
\begin{equation}\label{eq:galort}
A_{h,\alpha} ( u - u_{h,\alpha}, v) - \alpha (\nabla u,\nabla v )_{\Omega_h \setminus \Omega}
  = 0 \qquad \forall v \in V_h
\end{equation}
where the second term on the right hand side accounts for the inconsistent finite cell 
stabilization.
\end{tcolorbox}

\begin{proof} We first recall that for $\alpha = 0$ the method is 
consistent, i.e. the exact solution $u$ to (\ref{eq:problem-strong-a}--\ref{eq:problem-strong-b}) satisfies
\begin{equation}
A_{h,0}(u,v) = L_h(v) \qquad \forall v \in V_h
\end{equation}
since the standard Nitsche method with $\gamma = 0$ and $\tau=0$ is consistent 
and the added least squares terms are residual terms:
\begin{alignat}{3}
\Delta u  +  f &= 0 &
\quad 
&\Longrightarrow &
\quad
\tau h^2(\Delta u, \Delta v)_{\mcT_{h,\partial\Omega}\cap \Omega} &= -\tau h^2 (f, \Delta v)_{\mcT_{h,\partial\Omega} \cap \Omega} 
\\
u|_{\partial \Omega} &= g &
\quad 
&\Longrightarrow &
\quad 
2 h (\nabla_T u, \nabla_T v)_{\partial \Omega} &= 2 h (\nabla_T g,\nabla_T v)_{\partial \Omega} 
\end{alignat}
Next note that for $v\in V_h$ we have
\begin{align}
0 &= L_h(v) - A_{h,\alpha}(u_{h,\alpha}, v )
\\
&= A_{h,0}(u,v) - A_{h,\alpha}(u_{h,\alpha},v) 
\\
&= \bigl( A_{h,\alpha}(u,v) - \alpha(\nabla u, \nabla v)_{\Omega_h \setminus \Omega} \bigr)
- A_{h,\alpha}(u_{h,\alpha},v) 
\\
&=A_{h,\alpha}(u - u_{h,\alpha},v) - \alpha(\nabla u, \nabla v)_{\Omega_h \setminus \Omega}
\end{align}
\end{proof}

\subsection{Energy Norm}

Recall that $V = H^2(\Omega)$ and $V_h|_{\Omega} \subset V$, and define the 
following norms on $V$,
\begin{align}\label{eq:energy-norm} 
\tn v \tn^2_{h,\alpha} &= \tn v \tn^2_h + \alpha \| \nabla v \|^2_{\Omega_h \setminus \Omega}
\\ \label{eq:energy-norm-b} 
 \tn v \tn_{h}^2 &= \| v \|^2_{a_h} + \| v \|^2_{b_h}
 \\
\| v \|^2_{a_h} &= a_h(v,v) = \| \nabla v \|^2_\Omega + \tau h^2 \| \Delta v \|^2_{\mcT_{h,\partial\Omega} \cap \Omega}
\\
\| v \|^2_{b_h} &= b_h(v,v) = (2+\tau^{-1})h^{-1} \| v \|^2_{\partial \Omega} + 2 h \| \nabla_T v \|^2_{\partial \Omega}
\end{align}
\begin{rem}
Since $V_h$ may contain basis functions whose support has an arbitrarily small intersection with $\Omega$ 
the energy norm \eqref{eq:energy-norm-b}  is only guaranteed to be a semi-norm on $V_h$ in the case $\alpha=0$ while it is a norm on $V_h$ if $\alpha > 0$. For this reason we will need the finite cell stabilization to guarantee stability of the linear system of equations, see Section~\ref{section:linear-system}.
\end{rem}

\subsection{Results on LS-Nitsche}
Here we summarize the main results on the LS-Nitsche method \cite{ElfLarLar2019}, which is essentially our method \eqref{eq:method} in the case $\alpha=0$.
The central estimate in the analysis presented in \cite{ElfLarLar2019} is the following bound on the normal flux form, from which coercivity and continuity results are readily established, requiring only  minor adaptations to cover the present method when $\alpha > 0$.
\vspace{1em}

\begin{tcolorbox}[colback=white,boxrule=1.5pt]
\paragraph{Continuity of Normal Flux Form.} The following estimate holds 
\begin{align}\label{eq:trace}
\bigl| (\nabla_n v, w)_{\partial \Omega} \bigr|
&\lesssim 
\| v \|_{a_h} \| w \|_{b_h}\qquad v,w \in V
\end{align}
where the hidden constant is close to one for $h$ sufficiently small.

\tcblower

\paragraph{Coercivity.} For $\beta>0$ sufficiently large the form 
$A_{h,\alpha}$ is coercive 
\begin{equation}\label{eq:coercivity}
\tn v \tn_{h,\alpha}^2 \lesssim A_{h,\alpha}(v,v) \qquad v \in V
\end{equation}
 
\paragraph{Continuity.} The form $A_{h,\alpha}$ is continuous 
\begin{equation}\label{eq:continuity}
A_{h,\alpha}(v,w) \lesssim \tn v \tn_{h,\alpha}
\tn w \tn_{h,\alpha} \qquad v,w \in V 
\end{equation}

\end{tcolorbox}

\begin{rem}[Coercivity] Note that the coercivity holds on the full space $V = H^2(\Omega)$ and not only 
on the discrete space $V_h$ as is the case in the standard analysis of Nitsche's method. The 
reason is that the LS-Nitsche formulation circumvents the use of an inverse inequality in the 
proof of the coercivity estimate.
Also note that $\beta$ being sufficiently large for the LS-Nitsche method actually translates into a quite moderate number, for instance $\beta=5$ is sufficient. This is derived based on the good knowledge of the constant in \eqref{eq:trace}.
We refer to \cite{ElfLarLar2019} for details.
\end{rem}

\begin{rem}[Adaptation of Proofs] In \cite{ElfLarLar2019} the coercivity and continuity is established for $\alpha =0$ 
and it follows directly from the definition of $A_{h,\alpha}$ and the norm $\tn \cdot \tn_{h,\alpha}$ that these results also hold for $A_{h,\alpha}$. For instance, 
\begin{equation}
\tn v \tn^2_{h,\alpha} = \tn v \tn^2_{h,0} + \alpha \| \nabla v \|^2_{\Omega_h \setminus \Omega} 
\lesssim
A_{h,0}(v,v) + \alpha(\nabla v,\nabla v)_{\Omega_h\setminus \Omega} 
= 
A_{h,\alpha}(v,v) 
\end{equation}

\end{rem}

\subsection{Interpolation Error Estimates}
\label{section:interpolation}

Define the interpolant by 
\begin{equation}\label{eq:interpolant-def}
\pi_{h} : H^s(\Omega) \ni u \mapsto \Pi_{h} ( E u ) \in V_h  
\end{equation} 
where $\Pi_{h}: H^s(\Omega_h) \mapsto V_h$ is 
a spline space quasi-interpolant, see \cite{MR3202239}. We have the standard a priori error 
estimate 
\begin{equation}\label{eq:interpol-basic}
\| v - \pi_h v \|_{H^m(T)} \lesssim h^{s - m } \| v \|_{H^s(\mcN_h(T))}
\end{equation}
where $\mcN_h(T) = \cup_{\varphi \in B(T)} \supp(\varphi)$ and 
$B(T) = \{ \varphi \in B : T \subset \supp(\varphi)\}$.  We have the 
interpolation estimate 
\begin{equation}\label{eq:interpol-energy}
\tn v - \pi_h v \tn_{h,\alpha} \lesssim h^p \| v \|_{H^{p+1}(\Omega)}
\end{equation}
for $0\leq m \leq s \leq p+1$.

\subsection{Error Estimates}

\begin{tcolorbox}[colback=white,boxrule=1.5pt]
\begin{thm}[Energy Norm Error Estimate]
The following error estimate holds 
\begin{equation}\label{eq:error-est-energy}
\tn u - u_{h,\alpha} \tn_{h,\alpha} \lesssim h^{p} \| u \|_{H^{p+1}(\Omega)} 
+\alpha^{1/2} h^{1/2} \| u \|_{H^2(\Omega)}
\end{equation}
and with  
\begin{equation}\label{eq:epsilon-choice}
\alpha \lesssim  h^{2p-1} 
\end{equation}
we obtain
\begin{equation}\label{eq:error-est-energy-specific}
\tn u - u_{h,\alpha} \tn_{h,\alpha} \lesssim h^{p} \| u \|_{H^{p+1}(\Omega)} 
\end{equation}   
\end{thm}
\end{tcolorbox}

\begin{proof} We first have
\begin{align}
\tn u - u_{h,\alpha} \tn_{h,\alpha} 
&
\lesssim \tn u - \pi_h u \tn_{h,\alpha} 
+
\tn \pi_h u - u_{h,\alpha} \tn_{h,\alpha} 
\\
&
\lesssim h^p \| u \|_{H^{p+1} (\Omega)}
+
\tn \pi_h u - u_{h,\alpha} \tn_{h,\alpha} 
\end{align}
where we used the interpolation estimate (\ref{eq:interpol-energy}). In order to estimate the 
second term we start by using coercivity
\begin{align}
\tn \pi_h u - u_{h,\alpha} \tn_{h,\alpha} 
&\lesssim 
\sup_{v \in V_h} \frac{A_{h,\alpha}( \pi_h u - u_{h,\alpha}, v ) }{\tn v \tn_{h,\alpha}}
\\
&\lesssim \tn u - \pi_h u \tn_{h,\alpha} 
+
\sup_{v \in V_h} \frac{A_{h,\alpha}( u - u_{h,\alpha}, v ) }{\tn v \tn_{h,\alpha}}
\\
&\lesssim 
\tn u - \pi_h u \tn_{h,\alpha} 
+
\sup_{v \in V_h} \frac{\alpha(\nabla u, \nabla v)_{\Omega_h \setminus \Omega} }{\tn v \tn_{h,\alpha}}
\\
&\lesssim 
h^p \| u \|_{H^{p+1}(\Omega)}
+
\alpha^{1/2} h^{1/2} \| u \|_{H^2(\Omega)}
\end{align}
Here we added and subtracted an interpolant, used the identities 
\begin{equation}
A_{h,\alpha}( u - u_{h,\alpha}, v ) 
= 
A_{h,\alpha}( u , v ) - L_{h,\alpha} (v ) 
=
 A_{h,\alpha}( u , v ) - A_{h,0} (u,v )
=  
\alpha(\nabla u, \nabla v)_{\Omega_h \setminus \Omega}
\end{equation}
and finally the estimate 
\begin{equation}\label{eq:finite-cell-xx}
\alpha (\nabla u, \nabla v )_{\Omega_h \setminus \Omega} \lesssim \alpha^{1/2} h^{1/2} \| u \|_{H^2(\Omega)} \tn v \tn_{h,\alpha} 
\end{equation}
To verify (\ref{eq:finite-cell-xx}) we first note that 
\begin{equation}\label{eq:finite-cell-xx-a}
 \alpha (\nabla u, \nabla v )_{\Omega_h \setminus \Omega} 
\lesssim 
\alpha \|\nabla u \|_{\Omega_h \setminus \Omega} \|\nabla v \|_{\Omega_h \setminus \Omega} 
\lesssim 
\alpha^{1/2} \| \nabla u \|_{\Omega_h \setminus \Omega} \tn v \tn_{h,\alpha}
\end{equation}
By construction, and as illustrated in Figure~\ref{fig:geom-regions},
\begin{equation}
\Omega_h \setminus \Omega 
\subseteq U_h(\partial \Omega) \setminus \Omega
= 
\cup_{t \in [0,h]} \partial \Omega_t
\end{equation}
where $\partial \Omega_t = \{x \in U_h : \rho(x) = t \}$ and $\rho$ is the distance function 
to $\partial \Omega$. Using this notation we have the estimates 
\begin{align}\label{eq:finite-cell-xx-b}
\| \nabla u \|_{\Omega_h \setminus \Omega} 
&\lesssim 
\| \nabla u \|_{U_h(\partial\Omega)\setminus \Omega} 
\lesssim 
h^{1/2}   \sup_{t \in [0,h]} \|\nabla u \|_{\partial \Omega_t} 
\\ \label{eq:finite-cell-xx-c}
&\qquad 
\lesssim 
h^{1/2} \sup_{t \in [0,h]} \| u \|_{H^2(\Omega_t)} 
\lesssim 
 h^{1/2}  \| u \|_{H^2(\Omega_\delta)} 
\lesssim 
 h^{1/2}  \| u \|_{H^2(\Omega)}
\end{align}
where we used the trace inequality $\| v \|_{\partial \Omega_t} \lesssim \| v \|_{H^1(\Omega_t)}$ 
with $\Omega_t = \Omega \cup (\cup_{s\in [0,t]} \partial \Omega_s)$ and $v = \nabla u$, and at 
last we used the stability (\ref{eq:extension-continuity}) of the extension operator. Combining 
(\ref{eq:finite-cell-xx-a}) and (\ref{eq:finite-cell-xx-c}) we obtain (\ref{eq:finite-cell-xx}). Thus the proof is complete.
\end{proof}

\begin{tcolorbox}[colback=white,boxrule=1.5pt]
\begin{thm}[$L^2$ Error Estimate] The following estimate holds
\begin{equation}\label{eq:error-est-L2}
\| u - u_{h,\alpha} \|_{\Omega} \lesssim h^{p+1} \| u \|_{H^{p+1}(\Omega)} 
+  (\alpha^{1/2} h^{3/2} + \alpha h)  \| u \|_{H^2(\Omega)}
\end{equation}
and with $\alpha$ defined by (\ref{eq:epsilon-choice}) we obtain
\begin{equation}\label{eq:error-est-L2-specific}
\| u - u_{h,\alpha} \|_{\Omega} \lesssim h^{p+1} \| u \|_{H^{p+1}(\Omega)} 
\end{equation}
\end{thm}
\end{tcolorbox}

\begin{proof}
Let $\phi \in V$ be the solution to the dual problem
\begin{equation}
-\Delta \phi = \psi \quad \text{in $\Omega$}, \qquad \phi = 0 \quad \text{on $\partial \Omega$}
\end{equation}
Posing this problem in variational form, testing with $e=u-u_{h,\alpha}$, and using that $\phi_{\partial\Omega}=0$, gives us the identity
\begin{align} \label{eq:L2-decomp}
(e,\psi)_\Omega &= A_{h,0}(e,\phi) - \tau h^2 (\Delta e, \Delta \phi)_{\mcT_{h,\partial\Omega}\cap\Omega}
\end{align}
For the first term in \eqref{eq:L2-decomp} we recall that with $\alpha=0$ we obtain a consistent method and then proceeding with the obvious estimates
\begin{align}
A_{h,0}(e,\phi)
&= L_h(\phi) - A_{h,\alpha}(u_{h,\alpha},\phi )
+ \alpha (\nabla u_{h,\alpha},\nabla \phi )_{\Omega_h\setminus \Omega}
\\
&= L_h(\phi - \pi_h \phi ) - A_{h,\alpha}(u_{h,\alpha},\phi - \pi_h \phi)
+ \alpha (\nabla u_{h,\alpha},\nabla \phi )_{\Omega_h\setminus \Omega}
\\
&= A_{h,0}(u,\phi - \pi_h \phi ) - A_{h,\alpha}(u_{h,\alpha},\phi - \pi_h \phi)
+ \alpha (\nabla u_{h,\alpha},\nabla \phi )_{\Omega_h\setminus \Omega}
\\
&= A_{h,0}(u - u_{h,\alpha},\phi - \pi_h \phi ) 
+ \alpha (\nabla u,\nabla \pi_h \phi )_{\Omega_h\setminus \Omega}
\\
&\lesssim 
\tn u - u_{h,\alpha}\tn_{h,0} \tn \phi - \pi_h \phi\tn_{h,0} 
+ \alpha\|\nabla u\|_{\Omega_h \setminus \Omega} \|\nabla \pi_h \phi \|_{\Omega_h\setminus \Omega}
\\
&\lesssim 
h \tn u - u_{h,\alpha} \tn_{h,0} \| \phi \|_{H^2(\Omega)} 
+ \alpha h \| u \|_{H^2(\Omega)} \| \phi \|_{H^2(\Omega)}
\\
&\leq
( h \tn u - u_{h,\alpha} \tn_{h,\alpha}
+ \alpha h \| u \|_{H^2(\Omega)} ) \| \psi \|_\Omega
\end{align}
By the dual problem and the definition of the energy norm we readily get the bound
for the second term in \eqref{eq:L2-decomp}
\begin{align}
 -\tau h^2 (\Delta e, \Delta \phi)_{\mcT_{h,\partial\Omega}\cap\Omega}
 &\leq \tau h^2 \| \Delta e \|_{\mcT_{h,\partial\Omega}\cap\Omega}
 \| \Delta \phi \|_{\mcT_{h,\partial\Omega}\cap\Omega}
 \\&\lesssim
  \tau^{1/2}h
 \tn u - u_{h,\alpha} \tn_{h,\alpha}   \| \psi \|_{\Omega}
\end{align}
Finally, choosing $\psi=e=u-u_{h,\alpha}$ and using the energy norm estimate \eqref{eq:error-est-energy} we obtain 
\begin{align}
\| e \|_{\Omega} \lesssim 
h^{p+1} \| u \|_{H^{p+1}(\Omega)} + (\alpha^{1/2} h^{3/2} + \alpha h) \| u \|_{H^2(\Omega)}
\end{align}
which concludes the proof.
\end{proof}

\section{Stability of the Linear System of Equations} \label{section:linear-system}
We first recall some basic concepts and then we prove an inverse estimate and 
a Poincar\'e estimate, and finally a bound on the condition number of the stiffness 
matrix.
\subsection{Basic Definitions}
\begin{itemize}
\item
The stiffness matrix $\hatA$ is defined by
\begin{equation}\label{eq:stiffnessmatrix}
(\hatA \hatv,\hatw)_{\IR^N} = A_{h,\alpha}(v,w)\qquad \forall v,w \in V_h
\end{equation}
where $N$ is the dimension of $V_h$ and $\hatv \in \IR^N$ is the coefficient 
vector in the expansion 
\begin{equation}
v(x) = \sum_{i\in I} \hatv_i \varphi_i(x)
\end{equation}
of $v$ in terms of the basis functions $B = \{ \varphi_i : i \in I\}$. 
\item The following equivalence holds
\begin{equation}\label{eq:rn-eqv}
h^d \| \hatv \|^2_{\IR^N} \sim \| v \|^2_{\Omega_h}\qquad v \in V_h
\end{equation}

\item 
The condition number 
of $\hatA$ is defined by 
\begin{equation}\label{eq:cond-def}
\kappa=\mathrm{cond}(\hatA) = \frac{\lambda_{\max}}{\lambda_{\min}}
\end{equation}
where $\lambda_{\max}$ and $\lambda_{\min}$ are the largest and smallest eigenvalues of $\hatA$.
\end{itemize}

We shall now derive an estimate of the condition number. In addition to the coercivity 
(\ref{eq:coercivity}) and continuity (\ref{eq:continuity}) we will need an inverse estimate 
and a Poincar\'e estimate. 

\subsection{Inverse and Poincar\'e Estimates}
\begin{itemize}
\item The following inverse estimate holds
\begin{equation}\label{eq:cond-inverse}
\tn v \tn_{h,\alpha} \lesssim h^{-1} \| v \|_{\Omega_h}
\end{equation}
\begin{proof} Recall that 
\begin{equation}
\tn v \tn^2_{h,\alpha} = \| v \|^2_{a_h} + \alpha \| \nabla v \|^2_{\Omega_h \setminus \Omega} 
+ \| v \|^2_{b_h}
\end{equation}
Using standard inverse estimates 
\begin{align}
\| v \|^2_{a_h}  + \alpha \| \nabla v \|^2_{\Omega_h \setminus \Omega}
&=
\| \nabla v \|^2_{\Omega}  
+\tau  h^2 \| \Delta v \|^2_{\mcT_{h,\delta} \cap \Omega} 
+ \alpha \|\nabla v \|^2_{\Omega_h \setminus \Omega}
\\
&\leq (1 + \alpha) \| \nabla v \|^2_{\Omega_h} + \tau h^2 \|\Delta v \|^2_{\Omega_h}
\\
&\lesssim h^{-2} \| v \|^2_{\Omega_h} 
\end{align}
and using the inverse trace inequality $\| v \|^2_{\partial \Omega \cap T} 
\lesssim h^{-1} \| v \|^2_T$, $v \in V_h$, we get
\begin{align}
\| v \|^2_{b_h} 
&=
h^{-1} \| v \|^2_{\partial \Omega} + h \| \nabla_T v \|^2_{\partial \Omega}
\\
&\lesssim 
h^{-2} \| v \|^2_{\mcT_h(\partial \Omega)} + \| \nabla v \|^2_{\mcT_h(\partial \Omega)}
\\
&\lesssim 
h^{-2} \| v \|^2_{\mcT_h(\partial \Omega)} 
\\
&\lesssim 
h^{-2} \| v \|^2_{\Omega_h}
\end{align}
\end{proof}

\item The following Poincar\'e estimate holds 
\begin{equation}\label{eq:poincare}
\alpha h^{-2} \|  v \|^2_{\Omega_h \setminus \Omega}  + \| v \|^2_\Omega \lesssim \tn v \tn^2_{h,\alpha}
\end{equation}
and as a consequence 
\begin{equation}
\min(1,\alpha h^{-2} ) \| v \|^2_{\Omega_h} \lesssim \tn v \tn^2_{h,\alpha}
\end{equation}
\begin{proof} Using a standard Poincar\'e estimate we have 
\begin{equation}\label{eq:poincare-proof-a}
\| v \|_\Omega \lesssim \tn v \tn_{h,0}
\end{equation}
It remains to estimate $\alpha \| v \|^2_{\Omega_h \setminus \Omega}$, adding and subtracting 
the extension $v^e = v\circ p$, where $p:U_h \rightarrow \partial \Omega$ is the closest point 
mapping, and using the obvious bounds we obtain
\begin{align}
\alpha \| v \|^2_{\Omega_h \setminus \Omega}
&\lesssim 
\alpha \| v - v^e \|^2_{\Omega_h  \setminus \Omega}  
+
\alpha \| v^e \|^2_{\Omega_h \setminus \Omega}
\\
 &\lesssim 
\alpha h^2 \| \nabla v  \|^2_{\Omega_h  \setminus \Omega}  
+
\alpha \| v^e \|^2_{U_\delta \setminus \Omega} 
\\
 &\lesssim 
\alpha h^2 \| \nabla v  \|^2_{\Omega_h  \setminus \Omega}  
+
\alpha h  \| v \|^2_{\partial \Omega} 
\\
 &\lesssim 
h^2 \alpha \| \nabla v  \|^2_{\Omega_h  \setminus \Omega}  
+
(\alpha h^2)  h^{-1} \| v \|^2_{\partial \Omega} 
\\
 &\lesssim 
h^2 \max(1 , \alpha ) \tn v \tn^2_{h,\alpha}
\end{align}
which for $\alpha<1$ gives
\begin{equation}\label{eq:poincare-proof-b}
\alpha h^{-2} \| v \|^2_{\Omega_h \setminus \Omega} \lesssim \tn v \tn^2_{h,\alpha}
\end{equation} 
Together (\ref{eq:poincare-proof-a}) and (\ref{eq:poincare-proof-b}) give (\ref{eq:poincare}).
\end{proof}
\end{itemize}

\subsection{Condition Number Estimate}

\begin{tcolorbox}[colback=white,boxrule=1.5pt]
\begin{thm}[Condition Number Scaling] The condition number satisfies 
\begin{equation}
\kappa \lesssim \max(h^{-2}, \alpha^{-1})
\end{equation}
and with $\alpha \sim h^{2p-1}$, 
\begin{equation} \label{eq:cond-number-special}
\kappa \lesssim h^{-(2p-1)}
\end{equation}
\end{thm}
\end{tcolorbox}

\begin{proof} In view of the definition (\ref{eq:cond-def}) of the condition number $\kappa=\mathrm{cond}(\hatA)$
we need to estimate 
$\lambda_{\max}$ and $\lambda_{\min}^{-1}$.

\noindent 
{{\bf 1.} To estimate ${\lambda}_{{\max}}$} we 
use the definition (\ref{eq:stiffnessmatrix}) of the stiffness matrix to pass 
over to the bilinear form and then we use continuity (\ref{eq:continuity}) and 
the equivalence (\ref{eq:rn-eqv}) as follows
\begin{align}\label{eq:lambda-max}
\lambda_{\max} 
= \max_{v \in V_h \setminus 0} \frac{(\hatA \hatv, \hatv)_{\IR^N}}{\| \hatv \|^2_{\IR^N}}
= \max_{v \in V_h \setminus 0} \frac{A_{h,\alpha}(v,v)}{\| \hatv \|^2_{\IR^N}}
\lesssim 
 \max_{v \in V_h \setminus 0} \frac{\tn v \tn^2_{h,\alpha}}{\| v \|^2_{\Omega_h}}
\lesssim h^{-2}
\end{align}

\noindent
{{\bf 2.} To estimate ${\lambda}_{\min}^{-1}$} we again use the definition (\ref{eq:stiffnessmatrix}) 
of the stiffness matrix to pass to the bilinear form, then we use coercivity (\ref{eq:coercivity}), 
and the equivalence (\ref{eq:rn-eqv}), 
\begin{align}\label{eq:lambda-min}
\lambda_{\min}^{-1} 
= \max_{v \in V_h \setminus 0} \frac{\| \hatv \|^2_{\IR^N}}{(\hatA \hatv, \hatv)_{\IR^N}}
= \max_{v \in V_h \setminus 0} \frac{\| \hatv \|^2_{\IR^N}}{A_{h,\alpha}(v,v)}
\lesssim 
 \max_{v \in V_h \setminus 0} \frac{\| v \|^2_{\Omega_h}}{\tn v \tn^2_{h,\alpha}}
\lesssim (\min(1,\alpha h^{-2}))^{-1}
\end{align}

\noindent
{{\bf 3.}} Combining the estimates (\ref{eq:lambda-max}) and 
(\ref{eq:lambda-min}) we obtain 
\begin{align}
\kappa =
\lambda_{\max} \lambda^{-1}_{\min} 
\lesssim
h^{-2} (\min(1,\alpha h^{-2}))^{-1} 
=
h^{-2} \max(1,\alpha^{-1} h^{2})
=
\max(h^{-2}, \alpha^{-1} ) 
\end{align}
which completes the proof.
\end{proof}

\begin{rem}[Suboptimal Scaling]
The standard bound for the stiffness matrix condition number scaling in a finite element method is $h^{-2}$ but we here sacrifice this when choosing $\alpha$ such that the method is of optimal order. If we instead would choose $\alpha \sim h^2$ we would achieve the standard condition number scaling but on the other hand sacrifice the optimal order of the error estimates.
\end{rem}

\section{Numerical Experiments}
In this section we illustrate the stability and accuracy of the method using a 2D problem. Note however, that both the method and analysis are applicable also to problems in 3D or higher dimensions.

\paragraph{Implementation.}
The finite cell method in 2D is implemented in MATLAB and the linear system of equations is solved using a direct solver (MATLAB's \verb+\+ operator).
In all experiments we use tensor product quadratic B-spline basis functions, i.e. $C^1$-splines, on a uniform background grid. The geometry of the domain $\Omega$ is described as a high resolution polygon.
When computing condition numbers the stiffness matrix $\hat{A}$ is first converted to a full matrix on which MATLAB's \verb+eig+ is applied.

\paragraph{Parameter Choices.}
For the penalty parameter $\beta$ we choose a fix value $\beta=5$, since according to the analysis above this results in a method that is guaranteed coercive in any cut situation.
We vary the stabilization parameter $\tau$, which trades weight between the least-squares bulk stabilization and the effective penalization of the boundary condition in the Nitsche penalty term, see Table~\ref{table:tau}.
For comparison in the results below we also include experiments using standard Nitsche boundary conditions in the finite cell method, i.e. the method \eqref{eq:method} but without the new least squares stabilization terms. For a in some sense fair comparison, we let the effective Nitsche penalty parameter be the same so the Nitsche penalty term in both cases takes the form
\begin{align}
\beta (2 + \tau^{-1}) h^{-1} (v,w)_{\partial\Omega}
\end{align}
but note that we take no measures to ensure that the standard method is actually coercive.

For the finite cell stabilization in parts of the elements outside the domain, in accordance with \eqref{eq:epsilon-choice}, we set the finite cell stabilization parameter to $\alpha=0.001 h^{2p-1}$ where $p=2$ for the quadratic B-spline basis functions. According to the analysis this should yield optimal order convergence rates of $\mathcal{O}(h^{p})$ and $\mathcal{O}(h^{p+1})$ in energy respectively $L^2$ norms, but a suboptimal condition number scaling of $\mathcal{O}(h^{-(2p-1)})$.

\begin{table}[h]
    \centering
       \caption{
\emph{Nitsche Penalty.}
The effective Nitsche penalty parameter for various values of the least-squares bulk stabilization parameter $\tau$ when $\beta=5$.
} \label{table:tau}
    \begin{tabular}{cc}
        \toprule
  Least-squares stabilization & Nitsche penalty  \\
    $\tau$ &  $\beta (2+\tau^{-1})$  \\
        \midrule
1  & 15  \\
0.1 & 60 \\
0.01 & 510 \\
0.001  & 5010 \\
        \bottomrule
    \end{tabular}
\end{table}

\begin{figure}\centering
\begin{subfigure}[t]{0.32\linewidth}\centering
\includegraphics[trim=50 50 50 50, clip, width=0.935\linewidth]{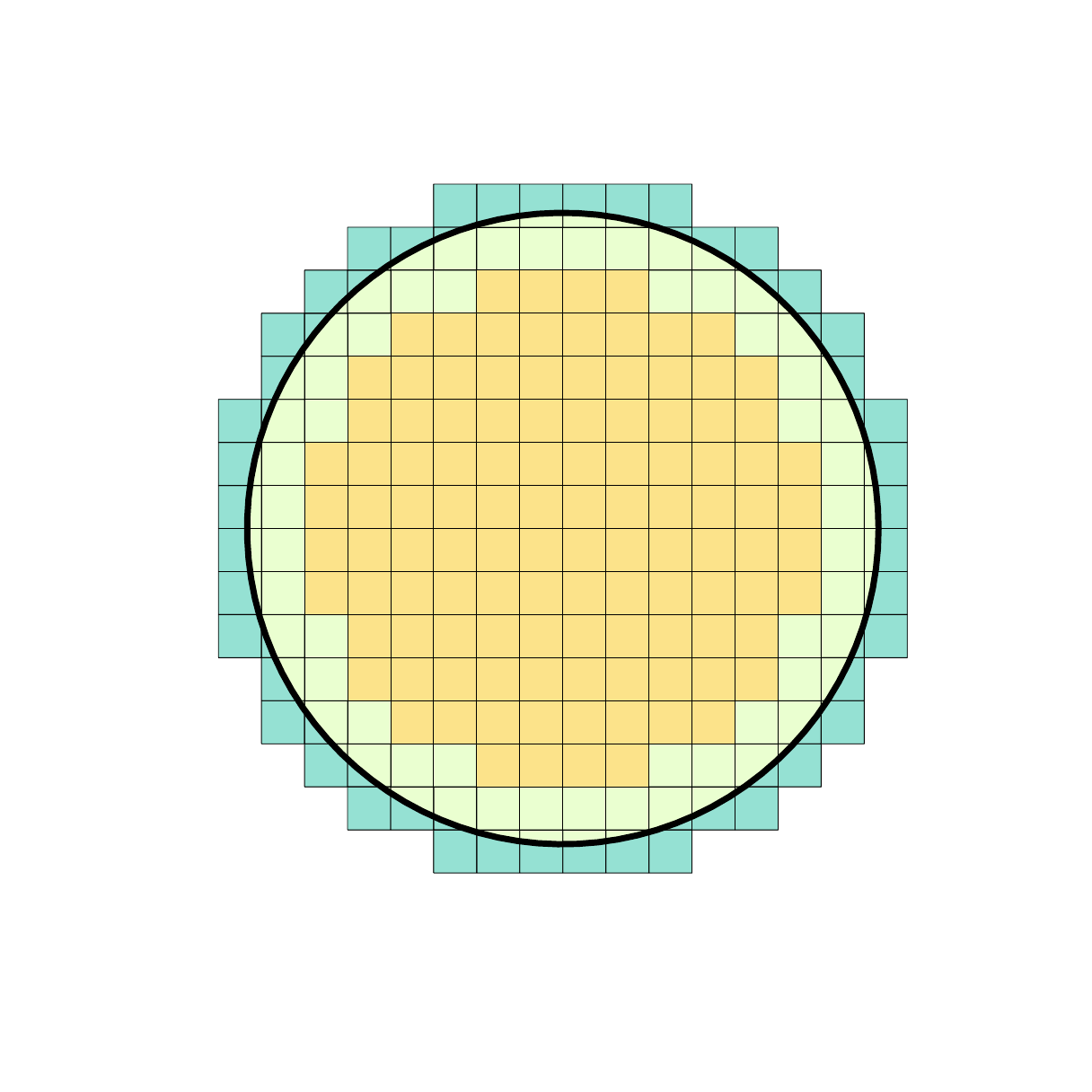}
\subcaption{$s=0$}
\end{subfigure}
\begin{subfigure}[t]{0.32\linewidth}\centering
\includegraphics[trim=50 50 50 50, clip, width=0.935\linewidth]{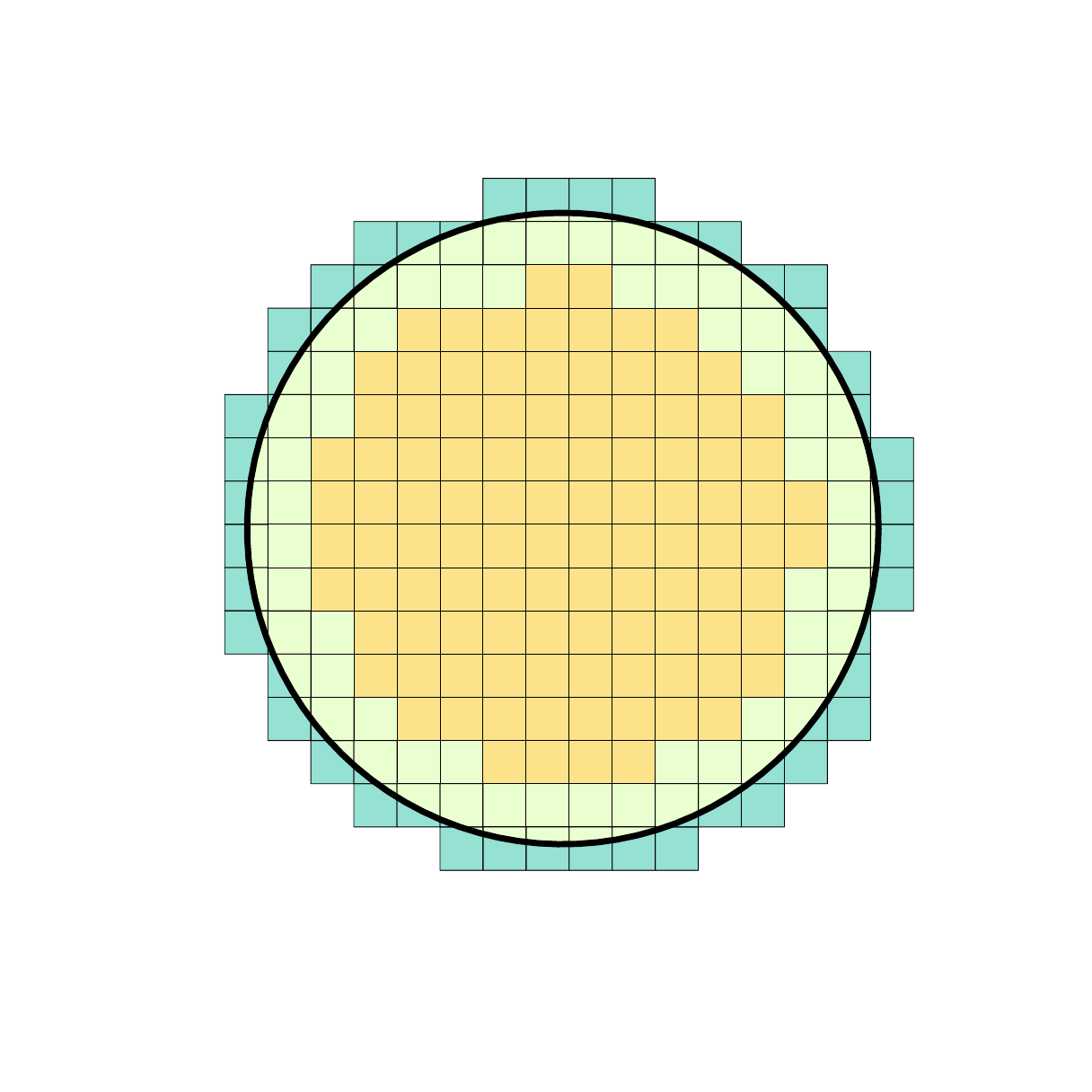}
\subcaption{$s=0.1$}
\end{subfigure}
\begin{subfigure}[t]{0.32\linewidth}\centering
\includegraphics[trim=50 50 50 50, clip, width=0.9\linewidth]{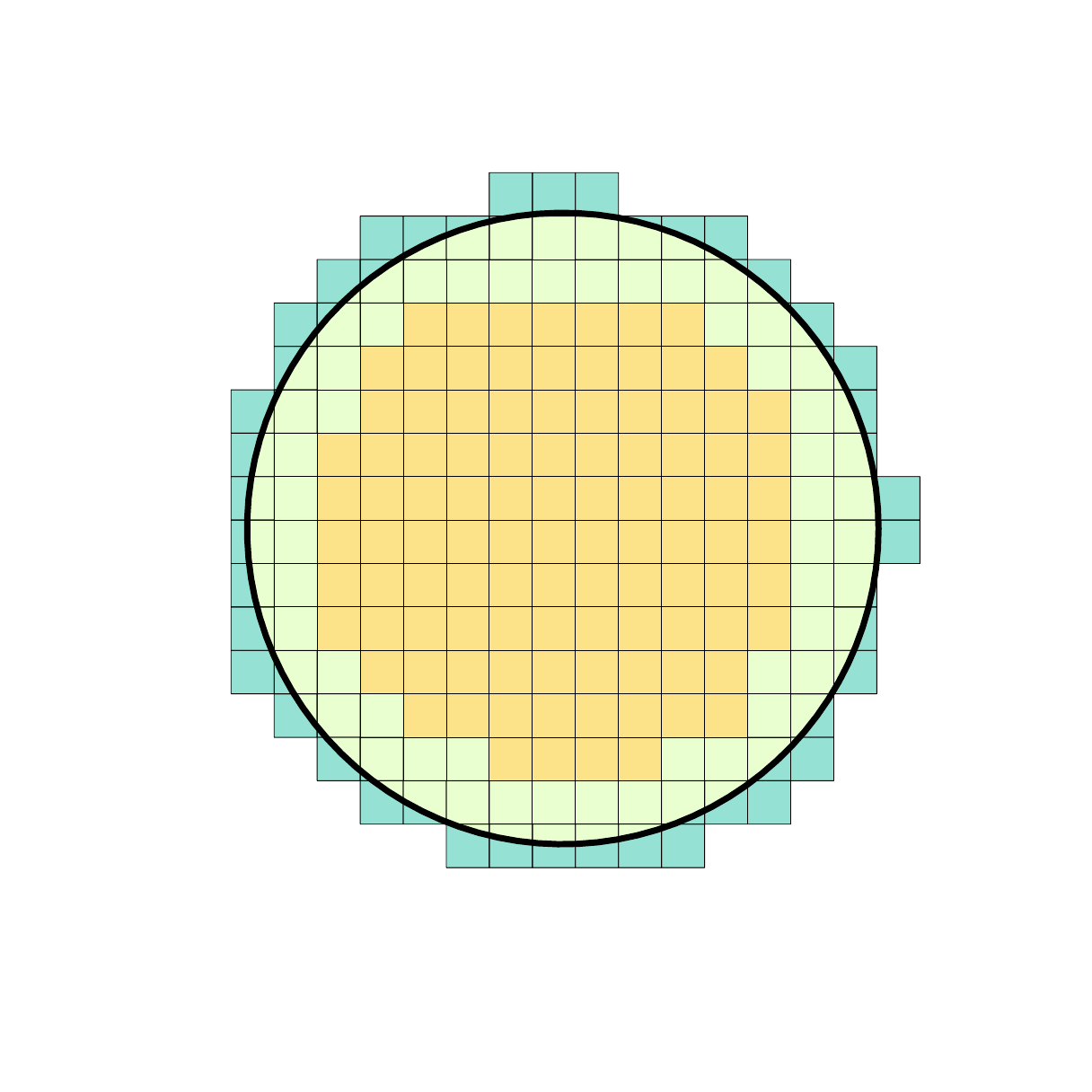}
\subcaption{$s=0.2$} \label{fig:disc-meshes:c}
\end{subfigure}
\caption{
\emph{Shifted Background Grid.}
Three example cut computational meshes ($h=0.1$) for the unit disc constructed by shifting the background grid $(sh,sh/3)$ where $s\in[0,1]$ is a parameter.
}
\label{fig:disc-meshes}
\end{figure}

\begin{figure}\centering
\begin{subfigure}[t]{0.40\linewidth}\centering
\includegraphics[trim=100 125 90 120, clip, width=0.85\linewidth]{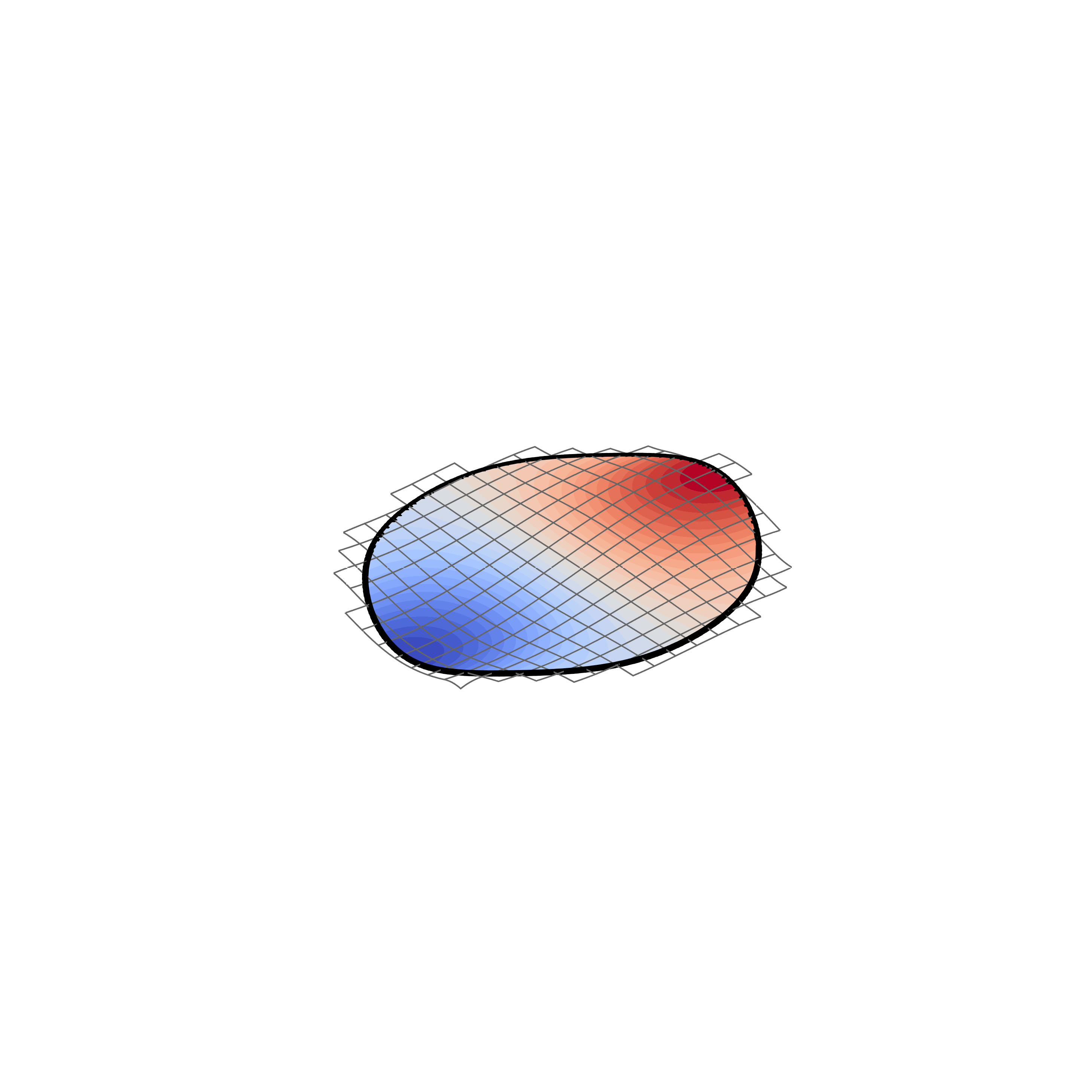}
\subcaption{LS-Nitsche}
\end{subfigure}
\begin{subfigure}[t]{0.40\linewidth}\centering
\includegraphics[trim=100 125 90 120, clip, width=0.85\linewidth]{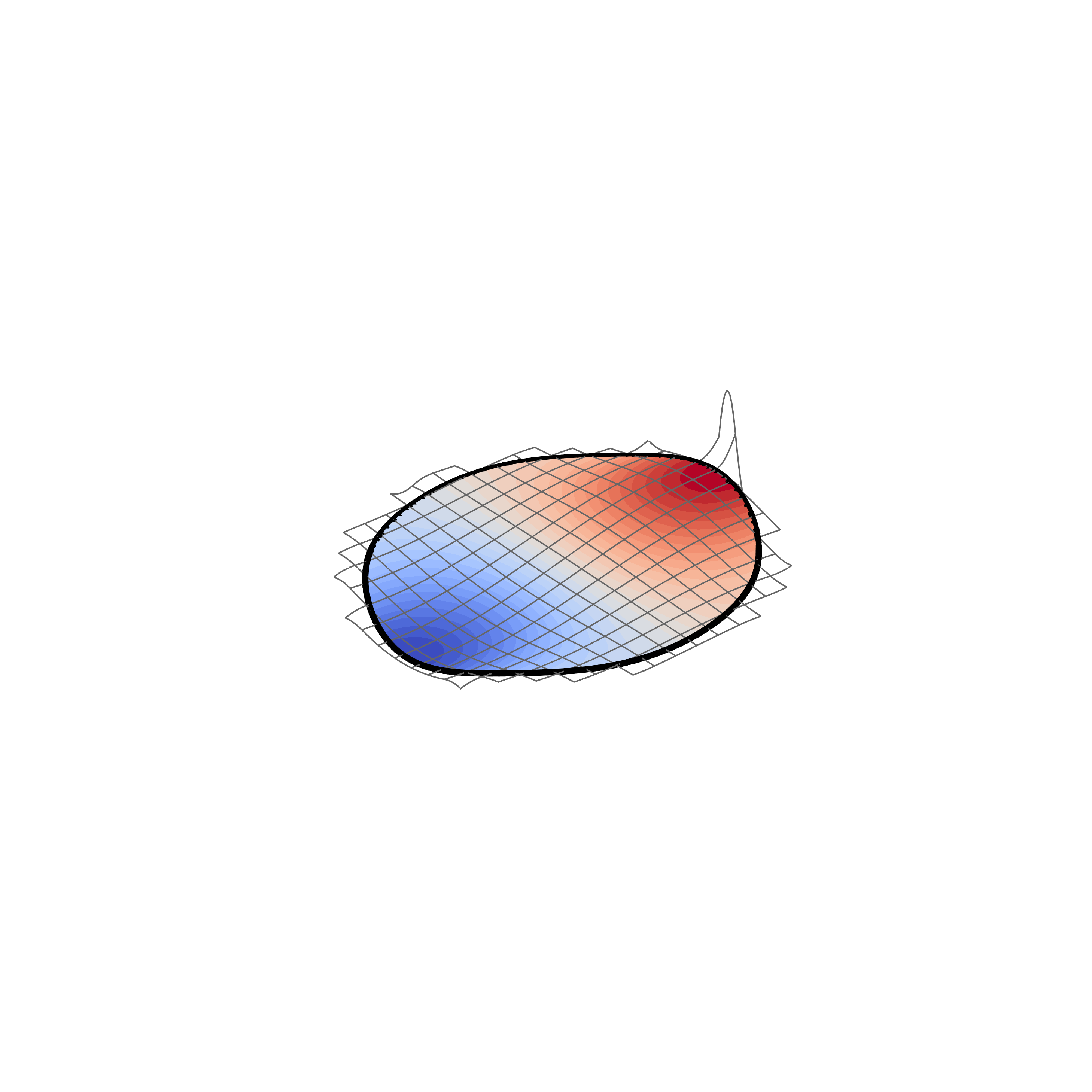}
\subcaption{Nitsche}
\end{subfigure}
\caption{
\emph{Solution.}
Numerical solutions for the manufactured problem on the computational mesh presented in
Figure~\ref{fig:disc-meshes:c}. Parameters here are $\beta=5$ and $\tau=0.1$ and the solution in (a) takes into account the new least-squares stabilization terms while the solution in (b) does not. The plotted solution outside the domain in (b) suggests some instabilities when not including the least-squares terms in this particular cut situation and choice of parameters.
}
\label{fig:example-solutions}
\end{figure}

\paragraph{Experiments.}
We manufacture a problem with known analytical solution on the unit disc $\Omega=\{(x,y): x^2+y^2 < 1\}$ via the ansatz
\begin{equation}
u(x,y) = \frac{1}{10}\left( \sin(2x) + x\cos(3y) \right)
\end{equation}
from which we derive data $f:\Omega\to\mathbb{R}$ and $g:\partial\Omega\to\mathbb{R}$.
To increase the probability of encountering problematic cut situations, for a given mesh size we perform computations on meshes extracted from 100 different shifts of the background grid and select the worst case. Specifically, the background grid is shifted $(sh,sh/3)$ where $s\in[0,1]$ is a parameter that we vary from 0 to 1 in 100 increments. Three example meshes are presented in Figure~\ref{fig:disc-meshes}. Finite cell solutions to the manufactured problem in one cut situation are presented in Figure~\ref{fig:example-solutions} with least-squares Nitsche stabilized boundary conditions respectively standard Nitsche boundary conditions. In this example we note an additional smoothness of the solution when we include the least-squares terms.

\begin{figure}\centering
\begin{subfigure}[t]{0.40\linewidth}\centering
\includegraphics[width=0.85\linewidth]{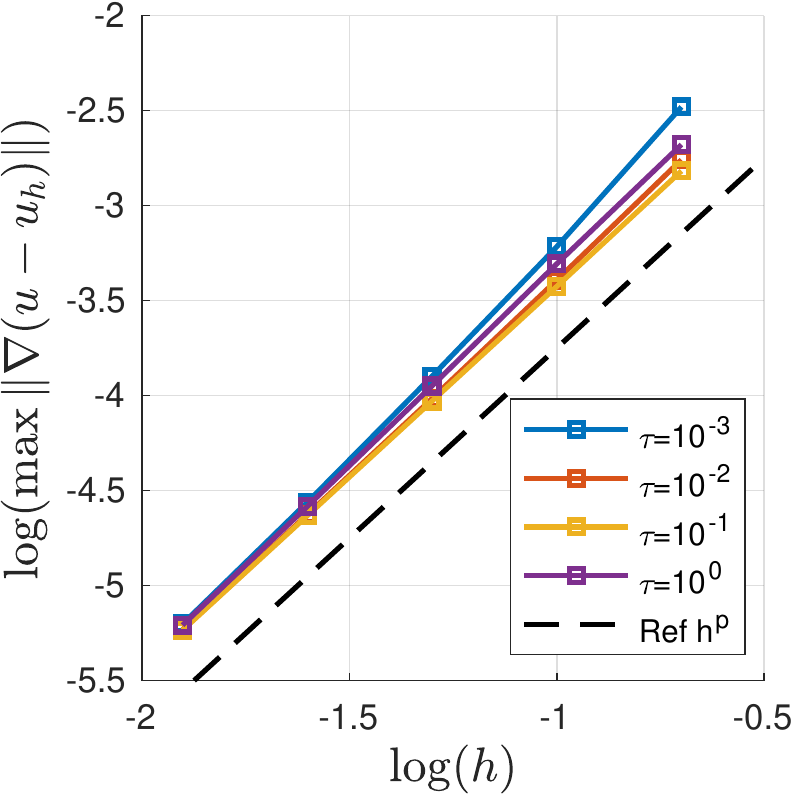}
\subcaption{LS-Nitsche}
\end{subfigure}
\begin{subfigure}[t]{0.40\linewidth}\centering
\includegraphics[width=0.85\linewidth]{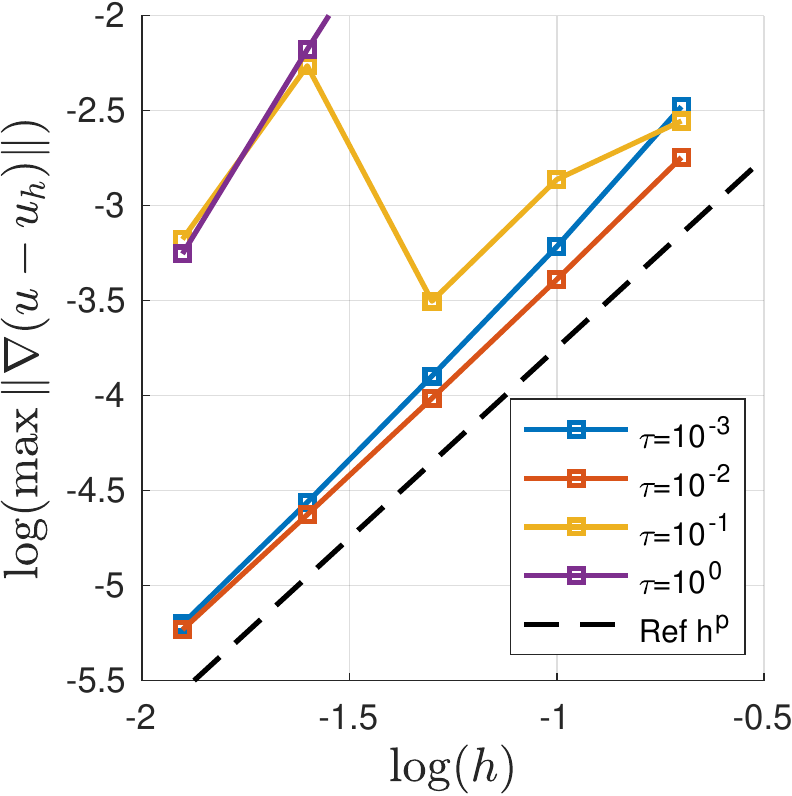}
\subcaption{Nitsche}
\end{subfigure}

\vspace{0.5em}
\begin{subfigure}[t]{0.40\linewidth}\centering
\includegraphics[width=0.85\linewidth]{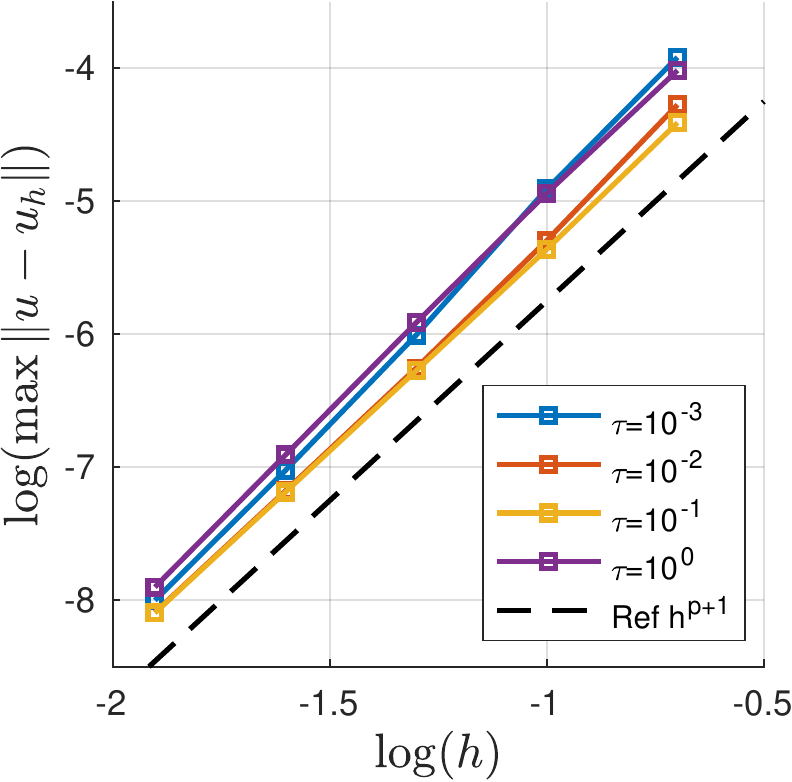}
\subcaption{LS-Nitsche}
\end{subfigure}
\begin{subfigure}[t]{0.40\linewidth}\centering
\includegraphics[width=0.85\linewidth]{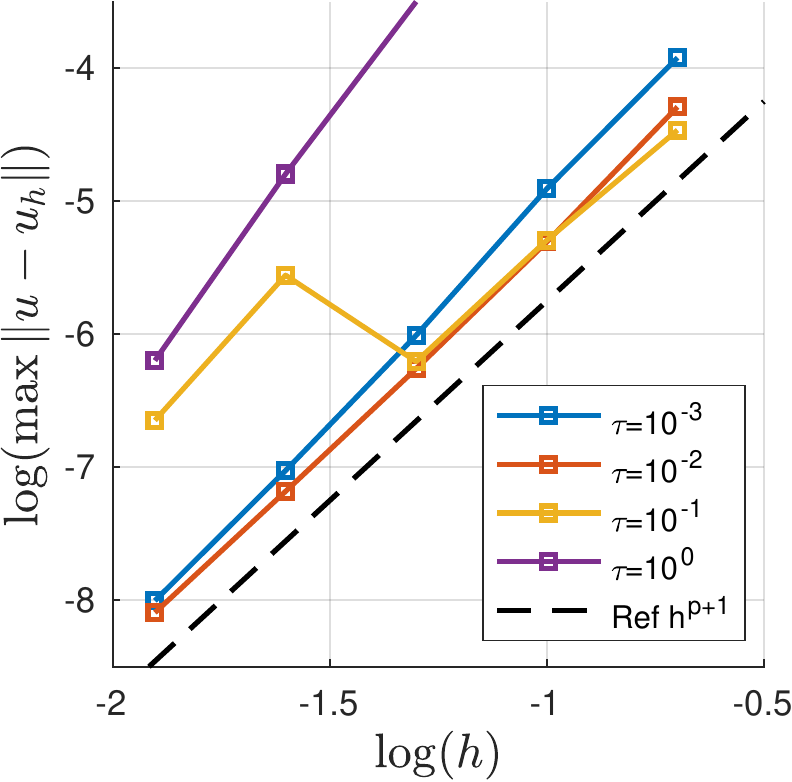}
\subcaption{Nitsche}
\end{subfigure}
\caption{
\emph{Convergence.}
Worst case convergence for a finite cell method with least squares stabilized Nitsche boundary conditions respectively standard Nitsche boundary conditions, both with the same effective penalty parameter.
}
\label{fig:worst-case-convergence}
\end{figure}

\emph{Convergence.} To illustrate the optimal order error estimates \eqref{eq:error-est-energy-specific} and  \eqref{eq:error-est-L2-specific} in
Figure~\ref{fig:worst-case-convergence} we present convergence results in $H^1$ seminorm and in $L^2$ norm for the manufactured problem using various values of $\tau$ for the finite cell method with least-squares stabilized Nitsche boundary conditions respectively standard Nitsche boundary conditions.
For each mesh size $h$ the largest error among the 100 different cut situations is presented, yielding a sort of worst case convergence. 
When using least-squares stabilized Nitsche boundary conditions the finite cell method achieves optimal order convergence for all values of $\tau$, while the method with standard Nitsche boundary conditions fails in this task for the larger values $\tau=1$ and $\tau=0.1$. The reason is most likely that the standard Nitsche method is then non-coercive in the worst cut situations. For lower values of $\tau$, which means a larger effective Nitsche penalty, we don't see this failure when using standard Nitsche boundary conditions. We do, however, see an increasing magnitude for the errors as $\tau$ is lowered which we attribute to locking due to the geometry being curved inside elements in combination with an increasing magnitude of the effective Nitsche penalty.
Overall we achieve the best performance in this experiment using $\tau=0.1$ and least-squares stabilized Nitsche boundary conditions.

\begin{figure}\centering
\begin{subfigure}[t]{0.40\linewidth}\centering
\includegraphics[width=0.85\linewidth]{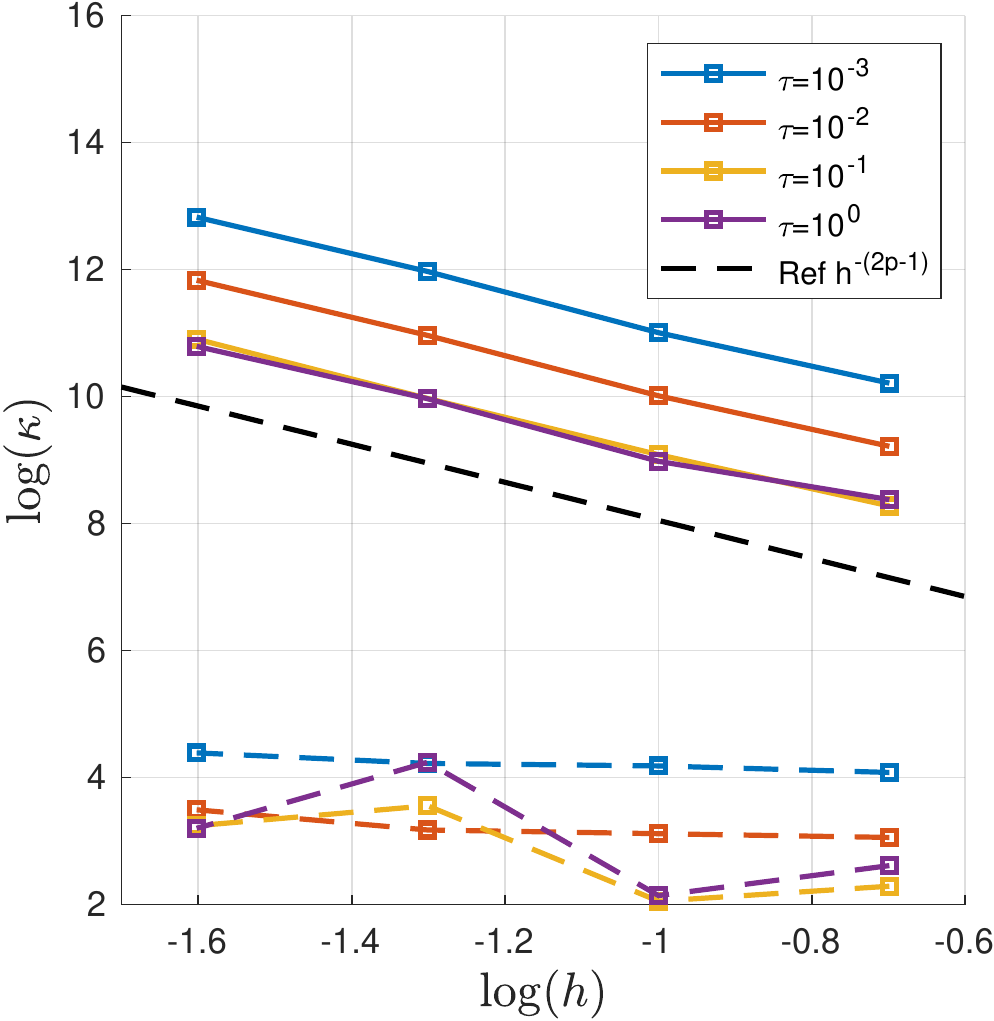}
\subcaption{LS-Nitsche}
\end{subfigure}
\begin{subfigure}[t]{0.40\linewidth}\centering
\includegraphics[width=0.85\linewidth]{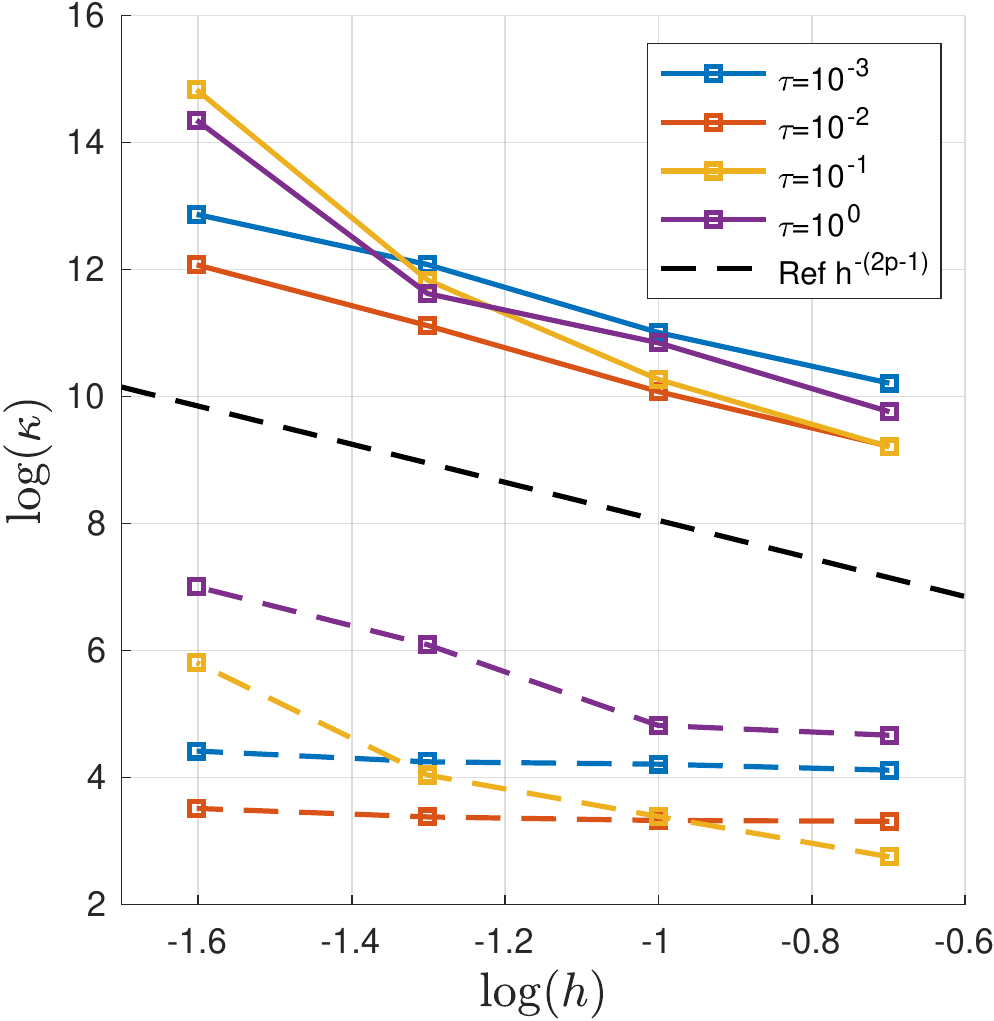}
\subcaption{Nitsche}
\end{subfigure}
\caption{
\emph{Condition Number Scaling.}
Worst case scaling of the condition number for a finite cell method with least squares stabilized Nitsche boundary conditions respectively standard Nitsche boundary conditions, both with the same effective penalty parameter. The dashed lines are the corresponding condition numbers after preconditioning using diagonal scaling.
}
\label{fig:worst-case-conditionnumber}
\end{figure}

\emph{Condition Number.}
To illustrate the condition number scaling \eqref{eq:cond-number-special}
we, analogously to the convergence experiment, present worst case scalings for the stiffness matrix condition number in Figure~\ref{fig:worst-case-conditionnumber}. In accordance with theory the finite cell method with least-square stabilized Nitsche boundary conditions produce condition numbers that scales as $\mathcal{O}(h^{-(2p-1)})$, while the method with standard Nitsche boundary conditions fails to achieve this bound when $\tau$ is large ($\tau=1,\tau=0.1$).
As in the case of the worst case convergence we suspect that the standard method is non-coercive in those situations. While we do not see such failures in the standard method when using smaller values of $\tau$, we note that choosing a smaller $\tau$ seems to have a significant effect of increasing the magnitude of the condition numbers. We also note that diagonal scaling seems to work very well as a preconditioner for the method.

\begin{figure}\centering
\begin{subfigure}[t]{0.32\linewidth}\centering
\includegraphics[width=0.935\linewidth]{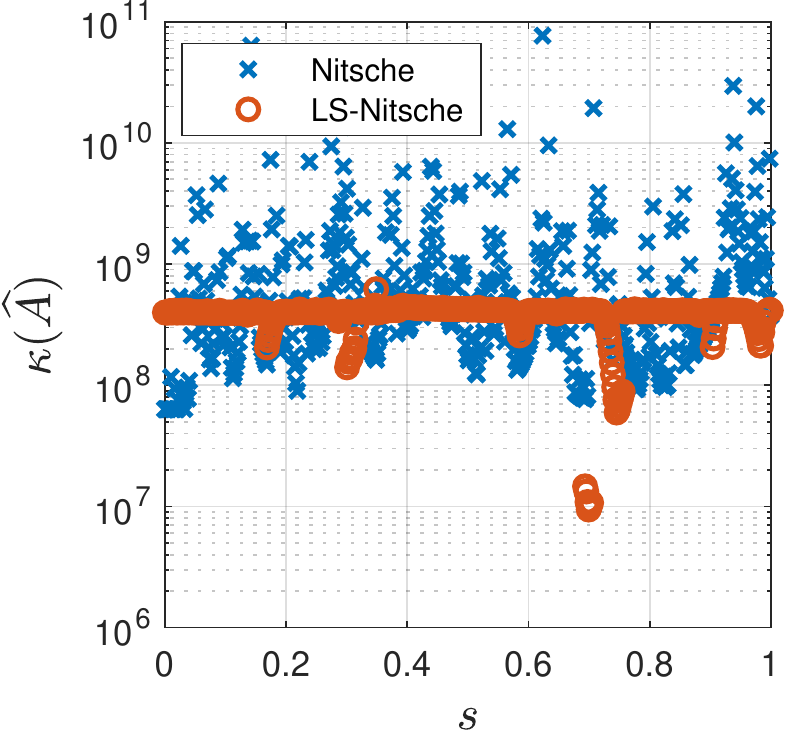}
\subcaption{$\tau=1$}
\end{subfigure}
\begin{subfigure}[t]{0.32\linewidth}\centering
\includegraphics[width=0.935\linewidth]{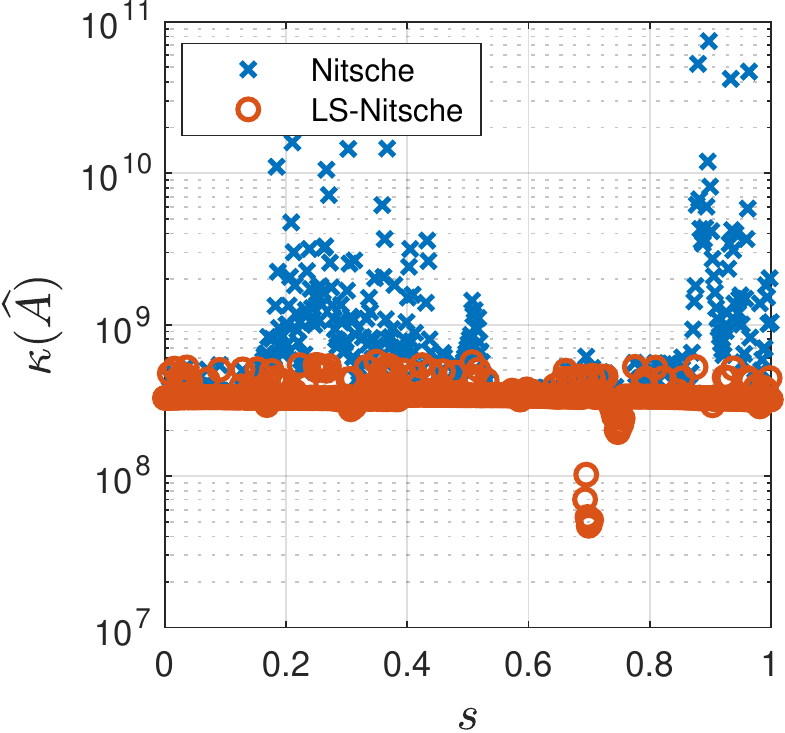}
\subcaption{$\tau=0.1$}
\end{subfigure}
\begin{subfigure}[t]{0.32\linewidth}\centering
\includegraphics[width=0.9\linewidth]{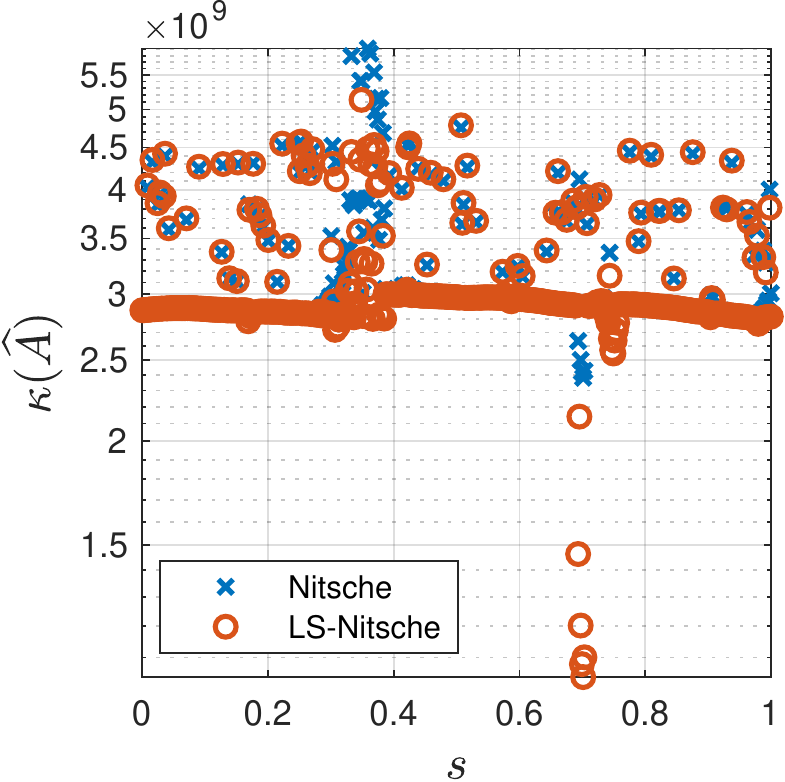}
\subcaption{$\tau=0.01$}
\end{subfigure}

\begin{subfigure}[t]{0.32\linewidth}\centering
\includegraphics[width=0.935\linewidth]{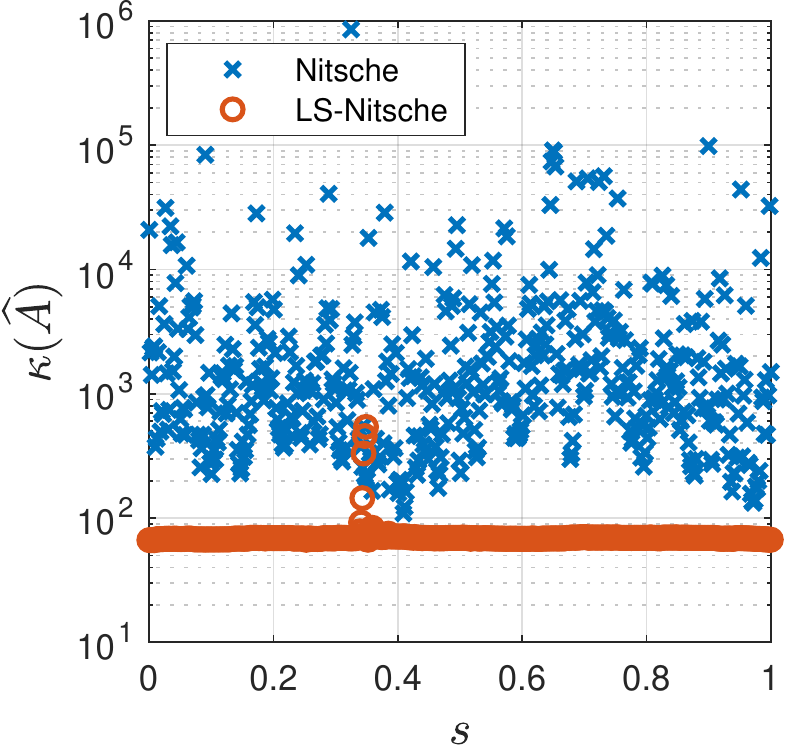}
\subcaption{$\tau=1$, preconditioned}
\end{subfigure}
\begin{subfigure}[t]{0.32\linewidth}\centering
\includegraphics[width=0.935\linewidth]{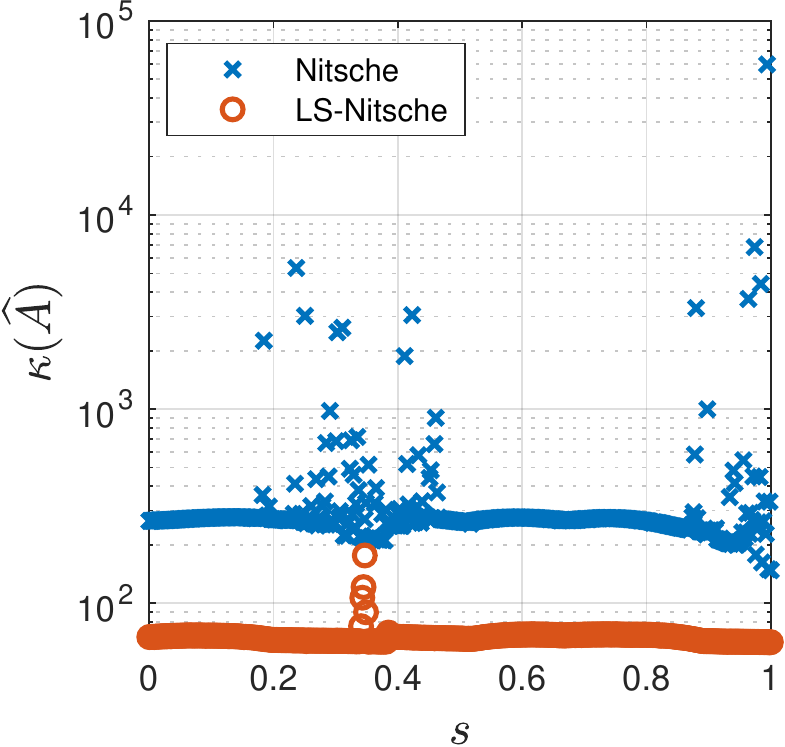}
\subcaption{$\tau=0.1$, preconditioned}
\end{subfigure}
\begin{subfigure}[t]{0.32\linewidth}\centering
\includegraphics[width=0.9\linewidth]{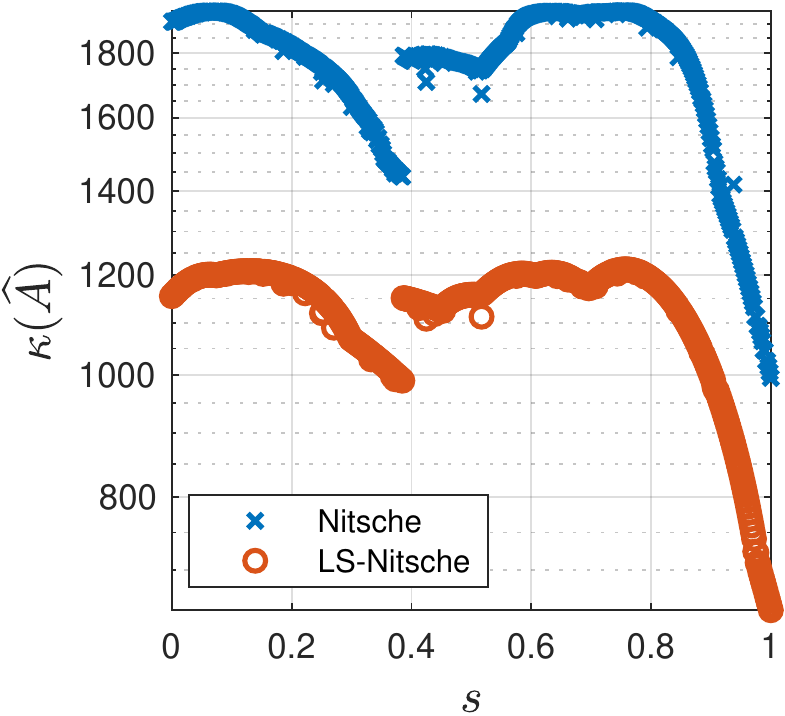}
\subcaption{$\tau=0.01$, preconditioned}
\end{subfigure}
\caption{
\emph{Condition Number Variation.}
Condition numbers for the stiffness matrices on the unit disc assembled using 500 different shifts of the background grid ($h=0.13$). 
}
\label{fig:conditionnumber-variation}
\end{figure}

In Figure~\ref{fig:conditionnumber-variation} we get a view of how the condition numbers vary over 500 different cut situations for the disc on a fixed mesh size, and we see the least-squares stabilization seems to have a very positive effect for larger values of $\tau$. Again, we see the effectiveness of diagonal scaling. In particular, the combination of finite cell stabilization, LS-Nitsche and diagonal scaling seems confidence inspiring. However, in (d)--(f) at around $s=0.35$, we see some minor disturbance and discontinuity, which could warrant further investigation.

\begin{figure}\centering
\begin{subfigure}[t]{0.32\linewidth}\centering
\includegraphics[width=0.935\linewidth]{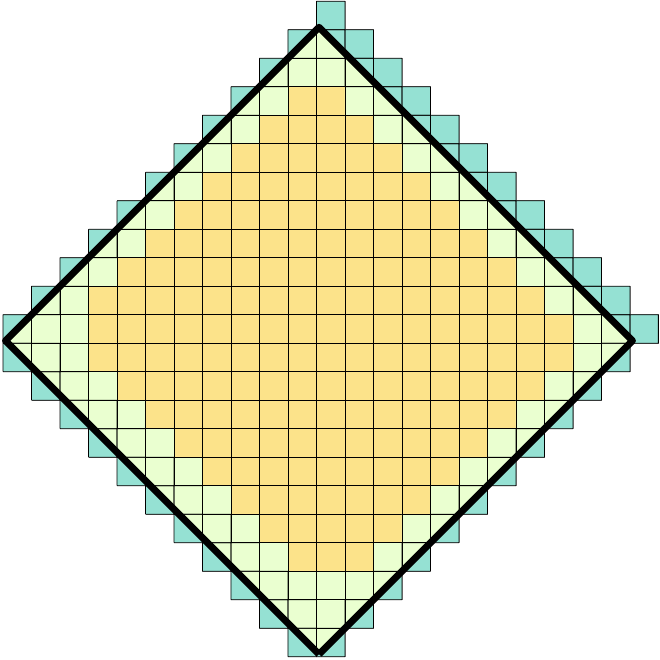}
\subcaption{$45^{\circ}$ unit square}\label{fig:spec-a}
\end{subfigure}
\begin{subfigure}[t]{0.32\linewidth}\centering
\includegraphics[width=0.935\linewidth]{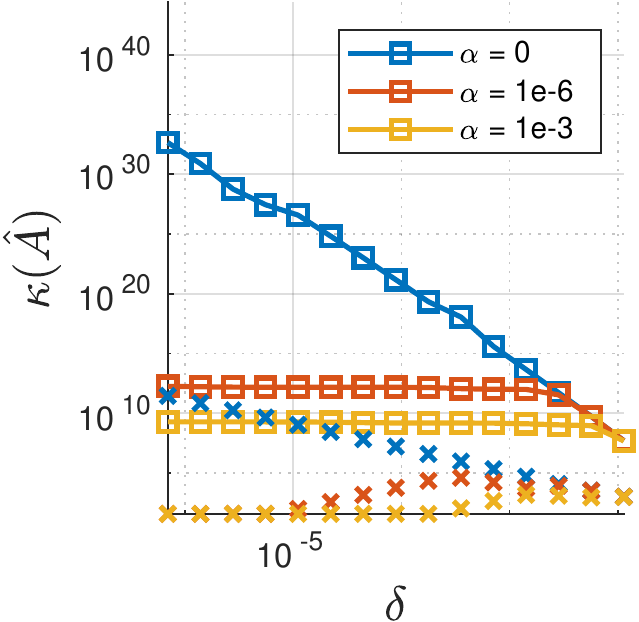}
\subcaption{LS-Nitsche}
\end{subfigure}
\begin{subfigure}[t]{0.32\linewidth}\centering
\includegraphics[width=0.935\linewidth]{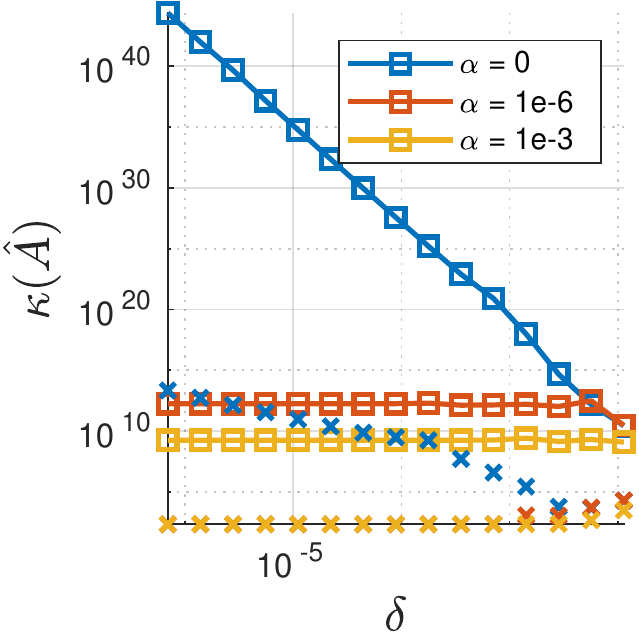}
\subcaption{Nitsche}\label{fig:spec-c}
\end{subfigure}
\\[1ex]

\begin{subfigure}[t]{0.32\linewidth}\centering
\includegraphics[trim=-25 -25 -25 -25, clip, width=0.935\linewidth]{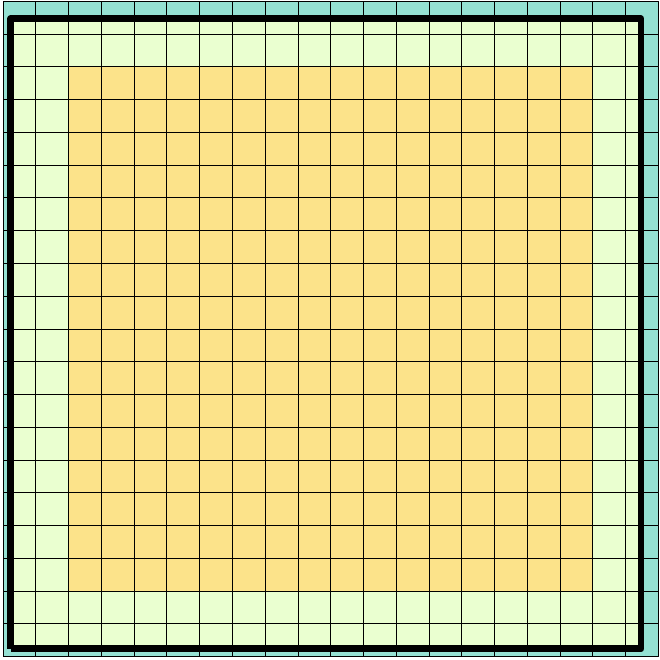}
\subcaption{Aligned unit square}
\end{subfigure}
\begin{subfigure}[t]{0.32\linewidth}\centering
\includegraphics[width=0.935\linewidth]{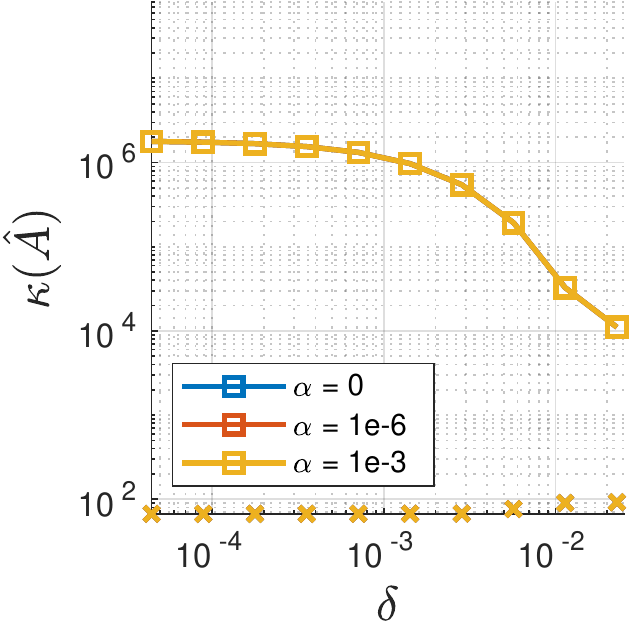}
\subcaption{LS-Nitsche}\label{fig:spec-d}
\end{subfigure}
\begin{subfigure}[t]{0.32\linewidth}\centering
\includegraphics[width=0.935\linewidth]{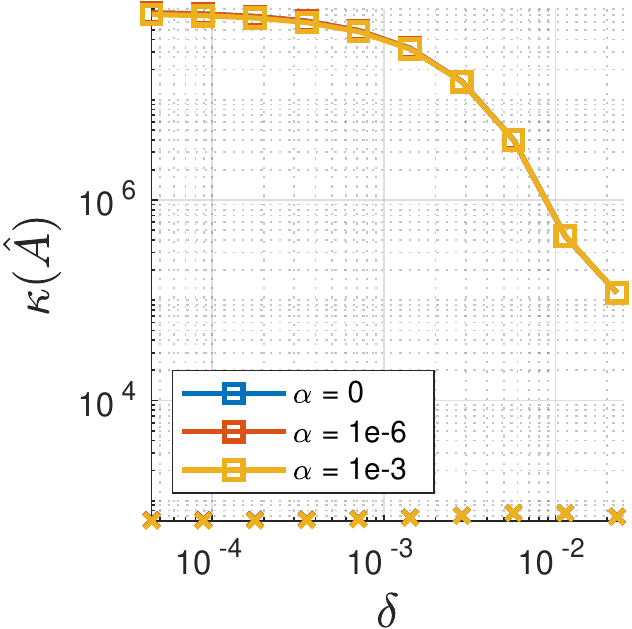}
\subcaption{Nitsche}\label{fig:spec-f}
\end{subfigure}
\caption{
\emph{Condition Numbers in Two Special Cases.}
Here two cases have been artificially constructed such that we create in a controlled way elements with an arbitrary small intersection with the domain depending on a parameter $\delta > 0$, where smaller $\delta$ means a smaller intersection.
The finite cell stabilization is here varied to show its effect on the condition numbers.
Note that to the values $\alpha = \{0,\text{1e-6},\text{1e-3}\}$ stated in the legends of these plots, the scaling $h^{2p-1}$ is also applied, and hence these values correspond to effective finite cell parameters of roughly $\alpha = \{ 0,10^{-10},10^{-7}\}$ in both cases. \emph{First case:} The geometry and mesh ($h=0.09$) are illustrated in (a), where the background grid is placed so that all cut elements are precisely half elements, after 
 which the geometry is shifted diagonally $(\delta h,\delta h)/\sqrt{2}$ creating arbitrarily small cut elements on the upper right side whose area is $\sim \delta^2$ of a full element. The condition numbers for this case are presented in (b) and (c), where crosses denote condition numbers after diagonal scaling.
\emph{Second case:} The geometry and mesh ($h=0.07$) are illustrated in (d), where the background grid is placed perfectly matching the geometry, after which the geometry is shifted diagonally $(\delta,\delta)$ creating arbitrarily small cut elements on the right and upper sides whose area is $\sim \delta$ of a full element and a corner element whose area is $\sim \delta^2$ of a full element. The condition numbers for this case are presented in (e) and (f), where crosses denotes the corresponding condition numbers after diagonal scaling.
}
\label{fig:conditionnumber-special}
\end{figure}

Finally, in Figure~\ref{fig:conditionnumber-special} we investigate the condition number in two special cases. Both cases are constructed such that the elements with the smallest intersection with the domain is $\sim \delta^2$ of a full element, and the parameter $\delta>0$ is varied. Further details of the two cases are given in the caption. In these experiments we also consider different values for the finite cell stabilization parameter, including zero, as well as preconditioning by diagonal scaling. For the first special case, Figures \ref{fig:spec-a}--\ref{fig:spec-c}, it becomes apparent that without finite cell stabilization, the condition number can become arbitrarily high, and this is not solved by the preconditioning. Including only a small amount of finite cell stabilization seems sufficient to bound the condition number. In the second special case, Figures \ref{fig:spec-d}--\ref{fig:spec-f}, the condition number of the system seems bounded even without finite cell stabilization. The effect of the least-squares stabilization terms does not seem to be essential in these experiments, but they generally give some improvement. Diagonal scaling on the other hand seems to work well as long as some finite cell stabilization is included.

\section{Conclusions}

In this contribution we have used least-squares stabilized Nitsche conditions to enforce Dirichlet boundary conditions in the finite cell method with $C^1$ splines. Our analysis and numerical experiments show the following:
\begin{itemize}
\item The method is symmetric and feature guaranteed coercivity in every cut situation, a property that does not depend on choosing a large penalty parameter.
 Since both the finite cell stabilization and the least-squares Nitsche stabilization terms are added element-wise, implementation in element based codes is straightforward.
\item In numerical experiments the method exhibits remarkable stability, both regarding convergence and conditioning, even in extreme cut
situations using a moderate size effective Nitsche penalty parameter ($\sim 10$).
\item Diagonal scaling seems to work quite well as a preconditioner for the method. While we during our  numerical experiments we have not seen a case where it breaks down when both the least-squares Nitsche stabilization terms and finite cell stabilization are included, we assume the method is not immune to the effect of cut basis functions being almost linearly dependent, in particular for higher order basis functions, see \cite{MR3610101,fritzphd}.
\end{itemize}
Overall, this combination of techniques for enforcing Dirichlet conditions in unfitted FEMs seems very promising.

\bibliographystyle{habbrv}
{\footnotesize{
\bibliography{ls-nitsche-refs}

\begin{thebibliography}{10}
\expandafter\ifx\csname url\endcsname\relax
  \def\url#1{\texttt{#1}}\fi
\expandafter\ifx\csname doi\endcsname\relax
  \def\doi#1{\burlalt{doi:#1}{http://dx.doi.org/#1}}\fi
\expandafter\ifx\csname urlprefix\endcsname\relax\def\urlprefix{URL }\fi
\expandafter\ifx\csname href\endcsname\relax
  \def\href#1#2{#2}\fi
\expandafter\ifx\csname burlalt\endcsname\relax
  \def\burlalt#1#2{\href{#2}{#1}}\fi

\bibitem{MR3202239}
L.~Beir\~{a}o~da Veiga, A.~Buffa, G.~Sangalli, and R.~V\'{a}zquez.
\newblock Mathematical analysis of variational isogeometric methods.
\newblock {\em Acta Numer.}, 23:157--287, 2014.
\newblock \doi{10.1017/S096249291400004X}.

\bibitem{Bur2010}
E.~Burman.
\newblock Ghost penalty.
\newblock {\em C. R. Math. Acad. Sci. Paris}, 348(21-22):1217--1220, 2010.
\newblock \doi{10.1016/j.crma.2010.10.006}.

\bibitem{BurClaHanLarMar2015}
E.~Burman, S.~Claus, P.~Hansbo, M.~G. Larson, and A.~Massing.
\newblock Cut{FEM}: discretizing geometry and partial differential equations.
\newblock {\em Internat. J. Numer. Methods Engrg.}, 104(7):472--501, 2015.
\newblock \doi{10.1002/nme.4823}.

\bibitem{DauDusRan2015}
M.~Dauge, A.~D\"{u}ster, and E.~Rank.
\newblock Theoretical and numerical investigation of the finite cell method.
\newblock {\em J. Sci. Comput.}, 65(3):1039--1064, 2015.
\newblock \doi{10.1007/s10915-015-9997-3}.

\bibitem{fritzphd}
F.~de~Prenter.
\newblock {\em Preconditioned iterative solution techniques for immersed finite
  element methods : with applications in immersed isogeometric analysis for
  solid and fluid mechanics}.
\newblock PhD thesis, Technische Universiteit Eindhoven, 2019.

\bibitem{MR3610101}
F.~de~Prenter, C.~V. Verhoosel, G.~J. van Zwieten, and E.~H. van Brummelen.
\newblock Condition number analysis and preconditioning of the finite cell
  method.
\newblock {\em Comput. Methods Appl. Mech. Engrg.}, 316:297--327, 2017.
\newblock \doi{10.1016/j.cma.2016.07.006}.

\bibitem{MR2458114}
A.~D\"{u}ster, J.~Parvizian, Z.~Yang, and E.~Rank.
\newblock The finite cell method for three-dimensional problems of solid
  mechanics.
\newblock {\em Comput. Methods Appl. Mech. Engrg.}, 197(45-48):3768--3782,
  2008.
\newblock \doi{10.1016/j.cma.2008.02.036}.

\bibitem{ElfLarLar2018}
D.~Elfverson, M.~G. Larson, and K.~Larsson.
\newblock Cut{IGA} with basis function removal.
\newblock {\em Adv Model Simul Eng Sci}, 5(6):1--19, 2018.
\newblock \doi{10.1186/s40323-018-0099-2}.

\bibitem{ElfLarLar2019}
D.~Elfverson, M.~G. Larson, and K.~Larsson.
\newblock A new least squares stabilized {N}itsche method for cut isogeometric
  analysis.
\newblock {\em Comput. Methods Appl. Mech. Engrg.}, 349:1--16, 2019.
\newblock \doi{10.1016/j.cma.2019.02.011}.

\bibitem{Fol1995}
G.~B. Folland.
\newblock {\em Introduction to partial differential equations}.
\newblock Princeton University Press, Princeton, 2:nd edition, 1995.

\bibitem{MR1941489}
A.~Hansbo and P.~Hansbo.
\newblock An unfitted finite element method, based on {N}itsche's method, for
  elliptic interface problems.
\newblock {\em Comput. Methods Appl. Mech. Engrg.}, 191(47-48):5537--5552,
  2002.
\newblock \doi{10.1016/S0045-7825(02)00524-8}.

\bibitem{MR2497337}
J.~Haslinger and Y.~Renard.
\newblock A new fictitious domain approach inspired by the extended finite
  element method.
\newblock {\em SIAM J. Numer. Anal.}, 47(2):1474--1499, 2009.
\newblock \doi{10.1137/070704435}.

\bibitem{JohLar2013}
A.~Johansson and M.~G. Larson.
\newblock A high order discontinuous {G}alerkin {N}itsche method for elliptic
  problems with fictitious boundary.
\newblock {\em Numer. Math.}, 123(4):607--628, 2013.
\newblock \doi{10.1007/s00211-012-0497-1}.

\bibitem{MR3310316}
D.~Kamensky, M.-C. Hsu, D.~Schillinger, J.~A. Evans, A.~Aggarwal, Y.~Bazilevs,
  M.~S. Sacks, and T.~J.~R. Hughes.
\newblock An immersogeometric variational framework for fluid-structure
  interaction: application to bioprosthetic heart valves.
\newblock {\em Comput. Methods Appl. Mech. Engrg.}, 284:1005--1053, 2015.
\newblock \doi{10.1016/j.cma.2014.10.040}.

\bibitem{MR3311666}
S.~Kollmannsberger, A.~\"{O}zcan, J.~Baiges, M.~Ruess, E.~Rank, and A.~Reali.
\newblock Parameter-free, weak imposition of {D}irichlet boundary conditions
  and coupling of trimmed and non-conforming patches.
\newblock {\em Internat. J. Numer. Methods Engrg.}, 101(9):670--699, 2015.
\newblock \doi{10.1002/nme.4817}.

\bibitem{Nit1971}
J.~Nitsche.
\newblock \"{U}ber ein {V}ariationsprinzip zur {L}\"{o}sung von
  {D}irichlet-{P}roblemen bei {V}erwendung von {T}eilr\"{a}umen, die keinen
  {R}andbedingungen unterworfen sind.
\newblock {\em Abh. Math. Sem. Univ. Hamburg}, 36:9--15, 1971.
\newblock \doi{10.1007/BF02995904}.

\bibitem{MR2377802}
J.~Parvizian, A.~D\"{u}ster, and E.~Rank.
\newblock Finite cell method: {$h$}- and {$p$}-extension for embedded domain
  problems in solid mechanics.
\newblock {\em Comput. Mech.}, 41(1):121--133, 2007.
\newblock \doi{10.1007/s00466-007-0173-y}.

\bibitem{MR2896903}
D.~Schillinger, A.~D\"{u}ster, and E.~Rank.
\newblock The {$hp$}-{$d$}-adaptive finite cell method for geometrically
  nonlinear problems of solid mechanics.
\newblock {\em Internat. J. Numer. Methods Engrg.}, 89(9):1171--1202, 2012.
\newblock \doi{10.1002/nme.3289}.

\bibitem{MR3358026}
D.~Schillinger and M.~Ruess.
\newblock The finite cell method: a review in the context of higher-order
  structural analysis of {CAD} and image-based geometric models.
\newblock {\em Arch. Comput. Methods Eng.}, 22(3):391--455, 2015.
\newblock \doi{10.1007/s11831-014-9115-y}.

\bibitem{MR3944486}
B.~Wassermann, S.~Kollmannsberger, S.~Yin, L.~Kudela, and E.~Rank.
\newblock Integrating {CAD} and numerical analysis: `dirty geometry' handling
  using the finite cell method.
\newblock {\em Comput. Methods Appl. Mech. Engrg.}, 351:808--835, 2019.
\newblock \doi{10.1016/j.cma.2019.04.017}.

\end{thebibliography}
}}

\bigskip
\paragraph{Acknowledgements.}
The first and last authors greatly acknowledge the funds provided by the Swedish Foundation
for Strategic Research Grant No.\ AM13-0029, the Swedish Research
Council Grants Nos.\  2013-4708, 2017-03911, and the Swedish
Research Programme Essence.
The second and third authors greatly acknowledge the funds provided by the Deutsche Forschungsgemeinschaft (DFG – German Research Foundation) Grants KO 4570/1-1 and RA 627/29-1.

\bigskip
\bigskip
\noindent
\footnotesize
\paragraph{Authors' addresses:}\ \smallskip

\smallskip
\noindent
Karl Larsson, \quad \hfill \addressumushort\\
{\tt karl.larsson@umu.se}

\smallskip
\noindent
Stefan Kollmannsberger,  \quad \hfill Chair of Computational Modeling and Simulation, TUM, Germany\\
{\tt stefan.kollmannsberger@tum.de}

\smallskip
\noindent
Ernst Rank, \quad \hfill Chair of Computational Modeling and Simulation, TUM, Germany\\
{\tt ernst.rank@tum.de}

\smallskip
\noindent
Mats G. Larson,  \quad \hfill \addressumushort\\
{\tt mats.larson@umu.se}

\end{document}